\documentclass[12pt,reqno]{article}

\begin{document}

\newcommand{\xc}[1]{}

\newcommand{\note}[1]{}

\newcommand{\innput}[1]{\input{#1}}

\newcommand{\pfbox}{QED}

\newcommand{\lab}[1]{\label{#1}\note{#1}}

\newcommand{\blr}{$\blacktriangleright\;$}

\newcommand{\inv}{{\bigcirc\!\!\!\!\!-1}}

\newcommand{\ba}{\begin{array}}
\newcommand{\ea}{\end{array}}

\newcommand{\f}[1]{\mbox{\boldmath $ #1$}}
\newcommand{\Bbb}[1]{\boldmath{#1}}

\newcommand \ruck{\hspace*{1mm}}
\newcommand \zuck{\hspace*{3cm}}

\newcommand{\N}{\Bbb{N}}
\newcommand{\Q}{\Bbb{Q}}
\newcommand{\st}{^{\bigcirc\!\!\!\!\!*}}

\newcommand{\bu}{\circ}
\newcommand{\sps}{\theta}
\newcommand{\bd}{\begin{quote}\begin{description}}
\newcommand{\ed}{\end{description}\end{quote}}
\newcommand{\ti}[1]{\tilde{#1}}
\newcommand{\ua}{\uparrow}
\newcommand{\da}{\downarrow}
\newcommand{\ey}{\emptyset}
 \newcommand{\io}{\iota}
\newcommand{\vep}{\varepsilon}\newcommand{\ep}{\epsilon}
\newcommand{\ov}{\overline}
\newcommand{\nperp}{\not \perp}
\newcommand{\cK}{\cal K}
\newcommand{\ded}{\vdash}
\newcommand{\dann}{\Rightarrow}
\newcommand{\et}{\wedge}
\newcommand{\imp}{\rightarrow  }
\newcommand{\eq}{\leftrightarrow}
 \newcommand{\ex}{\exists  }
 \newcommand{\all}{\forall  }

\newcommand{\Mod}{\sf Mod  }
\newcommand{\al  }{\alpha }
\newcommand{\be  }{\beta }
\newcommand{\ga  }{\gamma }
\newcommand{\de  }{\delta }
\newcommand{\Ga  }{\Gamma }
\newcommand{\De  }{\Delta }
\newcommand{\ka}{\kappa}
\newcommand{\la}{\lambda}
\newcommand{\si}{\sigma}
\newcommand{\Si}{\Sigma  }
\newcommand{\La}{\Lambda  }
\newcommand{\dale}{\vspace{.1in}}
\newcommand{\no}{\noindent}
\newcommand{\gdw}{\Leftrightarrow}
\newcommand{\s}{\subseteq}
\newcommand{\I}{\cal I}
\newcommand{\F}{\cal F}
\newcommand{\D}{\cal D}
\newcommand{\p}{^\perp}
\newcommand{\h}{^{\perp\!\perp}}
\newcommand{\om}{\omega}

\newtheorem{theorem}{Theorem}[section]
\newtheorem{proposition}[theorem]{Proposition}
\newtheorem{problem}[theorem]{Problem}
\newtheorem{exam}[theorem]{Example}

\newtheorem{lemma}[theorem]{Lemma}
\newtheorem{claim}[theorem]{Claim}
\newtheorem{corollary}[theorem]{Corollary}
\newtheorem{conjecture}[theorem]{Conjecture}
\newcommand{\ap}{\approx}
\newcommand{\app}{\!\!\approx\,\,}
\newcommand{\apr}{\!\!\approx_r\,\,}
\newcommand{\apl}{\!\!\approx_l\,\,}
\newcommand{\sub}{\subseteq}
\newcommand{\har}{{\cal L}(\hat{R}_{\hat{R}})}
\newcommand{\abs}{\dale\no}

\title{Proatomic  modular ortholattices: Representation
and equational theory}
\author{Christian Herrmann and  Michael S. Roddy}
\date{}

\maketitle
\begin{center} 
Presented on the Confernce\\ In memory of Herbert Gross\\
Locarno 1999\\
Note di Matematica e Fisica \textbf{10} (1999), 55--88

\end{center} 

\small
\noindent
{\it  \parbox{5cm}{
FB4 AG14 
TH Darmstadt
D64289 Darmstadt, Germany} \hfill
\parbox{7.5cm}{
Dept. of Mathematics, 
Brandon University\\
Brandon, Manitoba,
Canada 
R7A 6A9}}
\normalsize

\note{This is eq4.tex   of  Aug.22.99}

\begin{abstract}
We study  modular ortholattices in the variety  generated by
the finite dimensional ones from an equational and
geometric point of view. We relate this to coordinatization results.
\end{abstract}


\section{Introduction}
Modular lattices endowed with an orthocomplementation,
MOLs for short, 
 were  introduced by 
Birkhoff and von Neumann \cite{BVN:36}
 as  abstract
anisotropic orthogonal geometries. 
The cases of particular interest were the finite dimensional  
\cite{BVN:36} and  the continuous (von Neumann
\cite{NEU}) ones. These include the projection lattices of type $I_n$ resp.
type $II_1$ factors of von Neumann algebras.
According to Kaplansky \cite{KAP},
 completeness implies continuity and, in particular,
the absence of infinite families of pairwise perspective
and orthogonal elements ({\em finiteness}). This implies that in general
there is no completion. In particular, there is no obvious analogue
to ideal and filter lattices, the  basic tool in the
equational theory of lattices.

In our context, the most relevant result of that theory
is Frink's \cite{FRNK:46} embedding of a complemented
modular lattice in a subspace lattice of a projective space
and J\'{o}nsson's \cite{JON:54} supplement
that lattice identities are preserved under 
this construction. It easily follows that the lattice
variety generated by complemented modular lattices
is generated by its finite dimensional members (cf \cite{HH}).
 The r\^{o}le of finite resp.  finite dimensional MOLs
for the equational theory of MOLs was  discussed in 
Bruns \cite{GB:83}
and in Roddy  \cite{MR:87} focussing on a description of the 
lower part of the lattice of MOL-varieties.

In this paper our main objective are  the members
of the variety generated by finite dimensional MOLs. These will be called 
{\em proatomic} in view of the following
(where `geometric representation' refers to a projective
space with an anisotropic polarity). 
\begin{theorem}  The following are equivalent  for an MOL $L$
\begin{itemize}
\item[(1)] $L$ is proatomic 
\item[(2)] $L$ has an  atomic MOL-extension
\item[(3)] $L$ has a geometric representation
 \lab{atex}\end{itemize}
\end{theorem}
Our main tools are the MOL-construction method from 
Bruns and Roddy \cite{BR} and the concept of
orthoimplication from Herrmann and Roddy \cite{HR}.
The most prominent examples are the continuous geometries
constructed by von Neumann \cite{NEU2} from finite dimensional
inner product spaces. 
Also, 
we construct subdirectly irreducible proatomic MOLs
generated by an orthogonal 3-frame 
and of   arbitrarily large finite as well as infinite
height.

Quite a few questions remain unanswered - notably, whether
there is a non-proatomic MOL and whether every proatomic
MOL has an atomic extension within its variety.
Also, how to characterize $*$-regular rings with
a proatomic lattice of principal right ideals. These
and related questions are discussed in the final section.

As general references see \cite{B67,CD,GRO,HOL,KALM,KKW,MAE2,MAE3,NEU,SK,Wehr2}.
An excellent survey of complemented modular lattices
has been presented by Wehrung \cite{Wehr3}.
The most important concepts and results will be 
recalled in the sequel.

\section{Structure and coordinatization of MOLs}

\subsection{Complemented modular lattices}

All lattices will have smallest element $0$, treated as a constant.
Joins and meets will be written as $a+b$ and $ab$.
The {\em dimension} or {\em height} of a lattice  $L$ is the minimal length
of maximal chains where  length of a chain
is  cardinality   with one element deleted. 
$P_L$ denotes the set of atoms of $L$. $L$ is {\em atomic}
if for every $a>0$ there is $a \geq p \in P_L$.
And $L$ is {\em atomless} if it has no atoms, equivalently if
for all $a>0$ there is $a >b>0$.

Elements $a,b$ of a lattice form a {\em quotient} $a/b$ if $a \geq b$.
Then we have the interval sublattice 
$[a,b]=\{x \in L\mid a \leq x \leq b\}$ and we write 
$\dim[a,b] =\dim a/b$. The {\em height} of an element $a$ is
$\dim a/0$.
$a/b$  {\em transposes} down to  $c/d$ and $c/d$ up to $a/b$ if
$a=b+c$ and $d=bc$.  Quotients in the equivalence relation
 generated by transposed quotients are called {\em projective}
to each other.
Each lattice congruence is determined by its set of quotients
and closed under projectivity.

A lattice is {\em complemented} if it has bounds $0,1$ 
and if for every $a$ there is a {\em complement} $b$ such that
$ab=0$ and $a+b=1$. A lattice is {\em relatively complemented}
if each of its interval sublattices is complemented.
Any modular complemented lattice is such.

Elements $a,b$ of a lattice are {\em perspective}, $a \sim b$,
via $c$
if $c$ is a common complement  of $a,b$ in $[0,a+b]$.
In a complemented modular lattice,  $a \sim b$ via $d$ in $[ab,a+b]$ iff
$a \sim b$ via $c$ where $c$ is a complement of $ab$ in $[0,d]$
resp. $d =ab+c$.
Also,
 according to Lemma 1.4 in J\'{o}nsson \cite{J60},
if  $a \sim c \sim b$ and $a>0$
then there are $a \geq \ti{a} >0$ and $b \geq \ti{b}> 0$
such that $\ti{a} \sim \ti{b}$.

An ideal  is called {\em neutral} or 
a {\em p(erspectivity)-ideal}, if it is closed under
perspectivity.  
According to \cite{B67} p.78, for complemented modular lattices
   the neutral ideals $I$ are precisely the $0$-classes $I(\theta)$ of
lattice congruence relations $\theta$ - and determine those, uniquely:
 \[a/b \in \theta \mbox{ iff }  a=b+c \mbox{ for some } c \in I 
\;\;\mbox{ resp. }  ab' \in I 
\] 
Let $I(a)$ consist of all finite sums of $x_i$ perspective
to some $y_i \leq a$. By \cite{J60} Lemma 1.5 we have
that  $I(a)$ is the neutral ideal associated with the
congruence generated by $a/0$.
A lattice is {\em finitely subdirectly irreducible} if
the meet of any two nontrivial congruences is nontrivial.

\begin{proposition} \lab{neu}
Let $M$ be a subdirectly irreducible complemented modular lattice
 with minimal congruence $\mu$. Then \bd
\item 
$I(a) =I(\mu)$ for all $a/0 \in \mu$
\item
 $[0,b]$ is a simple lattice for each $b/0\in\mu$
\item 
For each $a>0$ there is $0<\tilde{a}\leq a$ with $\tilde{a}/0 \in \mu$.
\ed
A complemented modular lattice is finitely subdirectly irreducible iff
\bd
\item
For all  $a,b>0$ there
are $a \geq \tilde{a} >0$ and $b \geq \tilde{b} >0$
with $\tilde{a} \sim \tilde{b}$.
\ed
Every such is either atomic or atom-less.
\end{proposition}
{\em Proof}. Ad 1: $\mu$ is generated by any of its quotients.
Ad 2. Let $y<x\leq b$. Choose $a$ as complement of $y$ in $[0,x]$.
Since  $b \in I(\mu)=I(a)$ we have $b/0$ in the congruence of
$[0,b]$ generated by $a/0$, i.e. by $x/y$ cf Lemma 2.2 in \cite{J60}.\\
Ad 3. Let $c/d$ a generating quotient of $\mu$, w.l.o.g. $d=0$.
Then $c/0 \in con(a/0)$, i.e. $c/0$ has a proper sub-quotient 
projective to a sub-quotient $x/y$ of $a/0$.
But then $\tilde{a}/0 \in \mu$ with $a \geq \tilde{a}>0$
and a relative complement $\tilde{a}$ of $y$ in $[0,x]$.

Now, assume that $M$ is finitely subdirectly irreducible.
Then given $a,b>0$  we have $I(a) \cap I(b) \neq 0$ 
whence there $a\geq a_1 \sim c_1 \in I(b)$ and then 
$b \geq b_2 \sim c_2 \leq c_1$ and, by modularity,
$a_1 \geq a_2 \sim  c_2$ whence $\ti{a} \sim \ti{b}$ for some
$a_2\geq \ti{a}>0$ and $b_2\geq \ti{b}>0$. 
 If $L$ has an atom $a$, then each $b$ contains an atom perspective to $a$
in two steps. Conversely, we have $0<\ti{c} \in I(a) \cap I(b)$.
$\pfbox$

\subsection{Ortholattices}

An {\em ortholattice} is  a bounded  lattice, 
$L = (L; +, \cdot, 0, 1)$, together with an
orthocomplementation, i.e. a unary operation 
$\; ': L \mapsto L$ satisfying,
for all $x, y \in L$,
\[\mbox{
$x + x' = 1$,
$x \cdot x' = 0$, $x = x''$ and $x \leq y$ implies $y' \leq x'$.}\]  
Since 
the last
property, in the presence of the other ones,
 is equivalent to
 DeMorgan's laws
($(x + y)' = x' \cdot y'$ and its dual), 
this class of algebras forms a
variety, or equational class. {\em Modular}
  ortholattices  will be called
MOLs, for short. 
Examples are  Boolean algebras,   the height $2$ lattice MO$_\kappa$ with atoms
$a_\al,\,a'_\al\;(\al<\kappa)$ and
 orthocomplemented  non-desarguean planes,
e.g. arising by a free construction.

\xc{
\setlength{\unitlength}{0.00500000in}
\begingroup\makeatletter\ifx\SetFigFont\undefined%
\gdef\SetFigFont#1#2#3#4#5{%
  \reset@font\fontsize{#1}{#2pt}%
  \fontfamily{#3}\fontseries{#4}\fontshape{#5}%
  \selectfont}%
\fi\endgroup%
\begin{picture}(907,339)(0,-10)
\thicklines
\put(-27.000,93.000){\arc{304.631}{5.1173}{6.6881}}
\put(293.000,213.000){\arc{304.631}{1.9757}{3.5465}}
\put(93.000,133.000){\arc{304.631}{5.1173}{6.6881}}
\thinlines
\put(593,273){\blacken\ellipse{10}{10}}
\put(593,273){\ellipse{10}{10}}
\put(513,153){\blacken\ellipse{10}{10}}
\put(513,153){\ellipse{10}{10}}
\put(433,153){\blacken\ellipse{10}{10}}
\put(433,153){\ellipse{10}{10}}
\put(353,153){\blacken\ellipse{10}{10}}
\put(353,153){\ellipse{10}{10}}
\put(593,33){\blacken\ellipse{10}{10}}
\put(593,33){\ellipse{10}{10}}
\put(673,153){\blacken\ellipse{10}{10}}
\put(673,153){\ellipse{10}{10}}
\put(753,153){\blacken\ellipse{10}{10}}
\put(753,153){\ellipse{10}{10}}
\put(833,153){\blacken\ellipse{10}{10}}
\put(833,153){\ellipse{10}{10}}
\put(113,33){\blacken\ellipse{10}{10}}
\put(113,33){\ellipse{10}{10}}
\put(33,233){\blacken\ellipse{10}{10}}
\put(33,233){\ellipse{10}{10}}
\put(153,273){\blacken\ellipse{10}{10}}
\put(153,273){\ellipse{10}{10}}
\put(233,73){\blacken\ellipse{10}{10}}
\put(233,73){\ellipse{10}{10}}
\put(33,113){\blacken\ellipse{10}{10}}
\put(33,113){\ellipse{10}{10}}
\put(233,193){\blacken\ellipse{10}{10}}
\put(233,193){\ellipse{10}{10}}
\put(113,153){\blacken\ellipse{10}{10}}
\put(113,153){\ellipse{10}{10}}
\thicklines
\path(593,273)(513,153)
\path(513,153)(593,33)
\path(593,273)(673,153)
\path(673,153)(593,33)
\path(593,273)(433,153)
\path(433,153)(593,33)
\put(173.000,173.000){\arc{304.631}{1.9757}{3.5465}}
\path(593,273)(753,153)
\put(247,59){\makebox(0,0)[lb]{\smash{{{\SetFigFont{12}{14.4}{\familydefault}{\mddefault}{\updefault}$v'$}}}}}
\path(753,153)(593,33)
\path(593,273)(833,153)
\path(833,153)(593,33)
\path(593,273)(353,153)
\path(353,153)(593,33)
\path(33,233)(153,273)
\path(113,33)(233,73)
\path(233,193)(113,153)
\put(609,289){\makebox(0,0)[lb]{\smash{{{\SetFigFont{12}{14.4}{\familydefault}{\mddefault}{\updefault}$1$}}}}}
\put(608,2){\makebox(0,0)[lb]{\smash{{{\SetFigFont{12}{14.4}{\familydefault}{\mddefault}{\updefault}$0$}}}}}
\put(854,139){\makebox(0,0)[lb]{\smash{{{\SetFigFont{12}{14.4}{\familydefault}{\mddefault}{\updefault}$d'$}}}}}
\put(683,142){\makebox(0,0)[lb]{\smash{{{\SetFigFont{12}{14.4}{\familydefault}{\mddefault}{\updefault}$a'$}}}}}
\put(760,143){\makebox(0,0)[lb]{\smash{{{\SetFigFont{12}{14.4}{\familydefault}{\mddefault}{\updefault}$c'$}}}}}
\put(721,35){\makebox(0,0)[lb]{\smash{{{\SetFigFont{12}{14.4}{\familydefault}{\mddefault}{\updefault}$MO_3$}}}}}
\put(473,143){\makebox(0,0)[lb]{\smash{{{\SetFigFont{12}{14.4}{\familydefault}{\mddefault}{\updefault}$a$}}}}}
\put(398,143){\makebox(0,0)[lb]{\smash{{{\SetFigFont{12}{14.4}{\familydefault}{\mddefault}{\updefault}$c$}}}}}
\put(313,143){\makebox(0,0)[lb]{\smash{{{\SetFigFont{12}{14.4}{\familydefault}{\mddefault}{\updefault}$d$}}}}}
\put(100,0){\makebox(0,0)[lb]{\smash{{{\SetFigFont{12}{14.4}{\familydefault}{\mddefault}{\updefault}$0$}}}}}
\put(157,288){\makebox(0,0)[lb]{\smash{{{\SetFigFont{12}{14.4}{\familydefault}{\mddefault}{\updefault}$1$}}}}}
\put(0,231){\makebox(0,0)[lb]{\smash{{{\SetFigFont{12}{14.4}{\familydefault}{\mddefault}{\updefault}$v$}}}}}
\put(6,99){\makebox(0,0)[lb]{\smash{{{\SetFigFont{12}{14.4}{\familydefault}{\mddefault}{\updefault}$x$}}}}}
\put(241,202){\makebox(0,0)[lb]{\smash{{{\SetFigFont{12}{14.4}{\familydefault}{\mddefault}{\updefault}$x'$}}}}}
\put(59,139){\makebox(0,0)[lb]{\smash{{{\SetFigFont{12}{14.4}{\familydefault}{\mddefault}{\updefault}$x'v$}}}}}
\end{picture}
}

\abs 
{\em Orthomodular} lattices satisfy
only a special case of modularity:  $x=y+xy'$ for $y \leq x$.
It follows that $y \leq x$ generate a Boolean subalgebra and that 
 lattice congruences are
ortholattice congruences.
In particular, subdirect irreducibility depends only on the lattice structure
and we have  Prop.\ref{neu} for  MOLs, too.

\abs
Let $V(L)$ denote the ortholattice variety generated by $L$.
Any interval $[0,u]$ of an orthomodular lattice  is itself an orthomodular lattice
with complementation $x \mapsto x'u$ which is a homomorphic image of the subalgebra $[0,u] \cup [u',1]$ of $L$ whence in $V(L)$. 
Hence, by duality so are the intervals $[v,u]$.
We refer to these  as  {\em interval subalgebras}.
A {\em relative orthomodular lattice} is a lattice with an
orthomodular 
complementation on each of its interval sublattices,
such that each subinterval has the induced complementation. 
Thus, each orthomodular lattice $L$ can be considered as a relative one
and we have $M\in V(L)$ if and only if $M$
belongs to the relative variety of $L$.
In particular, an MOL, $L$, has the relative sub-MOL $L_{fin}$
which in turn can be considered as directed union
of the $[0,u]$, $u \in L_{fin}$.

\begin{lemma} Let $\sim $ be a reflexive binary relation  on
an orthomodular lattice $L$ which is   compatible with the lattice
operations (i.e.  a sublattice of $L^2$).
If $\sim $ is symmetric or compatible with the orthocomplement
(i.e. a subalgebra of $L^2$) then $\sim $  is a congruence relation
of $L$ \lab{tol}
\end{lemma}
Proof. If $\sim $ is also symmetric (i.e. a lattice {\em tolerance})
then we have $a \sim  b$ iff $a+b \sim  ab$.
Namely,
 $a+b \sim b+b =b$ and $a+b \sim a+a=a$ from $a \sim b$ resp.
$b \sim a$ whence   $a+b =(a+b)(a+b) \sim ab$.
Conversely, from $a+b\sim ab$ if follows 
 $a+b =a+b+b \sim ab+b=b$ and similarly $a+b \sim a$
whence $a \sim a+b$  and  $a=a(a+b) \sim (a+b)b=b$.
Therefore, from $a \sim b$ with $c=a+b \sim ab=d$ it follows
 $cd'\sim 0$ and, since
$d'= cd'+c'$ by orthomodularity, $a'+b'=d' \sim c'=a'b'$ whence $a' \sim b'$.

This means that $\sim$ is a subalgebra of $L^2$, in any case.
Now, recall that
 $p(x,y,z)=(x+((y+z)y')(z+((x+y)y')$ is a Mal'cev term 
for orthomodular lattices, i.e. $p(x,x,z)=z$ and $p(x,z,z)=x$.
Thus, according to  the Goursat-Lambek Lemma \cite{La} p.10
we have symmetry and  transitivity, too. Indeed, 
from $y \sim y$, $x \sim y$, and $ x \sim x$ it follows
$y =p(y,x,x)\sim p(y,y,x)=x$ and
from  $x\sim y$,
$y\sim y$, and $y\sim z$ it follows
$x=p(x,y,y)\sim   p(y,y,z)=z$. $\pfbox$

\begin{corollary} \lab{toll} A set ${\cal Q}$ of quotients in an
orthomodular lattices is the set of quotients of a congruence relation
(i.e. $a \theta b$ iff $(a+b)/(ab) \in {\cal Q}$) if and only if
it contains all $a/a$ and is closed under subquotients, transposes and
\[ a/c, b/c \in {\cal Q} \mbox{ implies } (a+b)/c \in {\cal Q}, \quad 
c/a,c/b \in {\cal Q} \mbox{ implies } c/(ab) \in {\cal Q}\]
\end{corollary}
{\em Proof}.  According to  \cite{Ban} 
$\theta$ is a lattice tolerance. Also,
the transitivity of ${\cal Q}$ is immediate from the
existence of relative complements. $\pfbox$
The most prominent example of a congruence on an MOL and its
neutral ideal are 
 \[ a\,\theta_{fin}\,b \;\Leftrightarrow\; \dim[\,ab,a+b]<\infty 
\quad \quad I=L_{fin}=\{a\in L\mid \dim[ 0,a]<\infty\}\]

\subsection{Review of coordinatization}

Let $n \geq 3$ fixed. An {\em $n$-frame}, in the sense of von Neumann \cite{NEU}, in a lattice  $L$ is a list $\f{a}:
a_{i}, a_{ij}, 1 \leq i,j \leq n, i \neq j$ of elements of $L$ such that for
any 3 distinct  $j,k,l$
\[ a_j \sum_{i \neq j}a_i =\prod_i a_i =a_ja_{jk},\;\;a_j+a_{jk}=a_j+a_k,\;\;
 a_{jl}=a_{lj}=(a_j+a_l)(a_{jk}+a_{kl})  .\]
The frame is {\em spanning} in $L$ if $\prod_i a_i =0_L$ and 
$\sum_i a_i =1_L$.
The {\em coordinate domains} associated with the frame $\f{a}$  are
 \[R_{ij}=R(L,\f{a})_{ij} = \{r \in L\,|\, ra_j=a_ia_j,\;\; r+a_j=a_i+a_j \} \quad i \neq j.\]
Now assume that $L$ is modular and $n \geq 4$ or in case $n=3$  
assume the Arguesian law of J\'{o}nsson \cite{JON:54}.
According to von Neumann \cite{NEU} and 
Day and Pickering \cite{DP},
using lattice polynomials $\oplus_{ij}$, $\ominus_{ij}$, $\otimes_{ij}$ in $\f{a}$, each of these
 can be turned into a ring  with zero $a_i$ and unit $a_{ij}$ such that
there are ring  isomorphism of
$R_{ij}$     onto $R_{ik}$ and $R_{kj}$ respectively
\[\pi_{ijk}r=r_{ik}= (r+a_{jk})(a_i+a_k), \quad
\pi_{jik}r=r_{kj}= (r+a_{ik})(a_k+a_j) .\]
Thus, we can speak of the {\em ring} $R(L,\f{a})$.
The operations on $R_{ij}$ can be defined with just one auxiliary
index $k$ and the result does not depend on the choice of $k$.
In particular,  the
multiplication on $R_{ik}$ is given by  
\[ (s\cdot r)_{ik} =    (r_{ij} +s_{jk})(a_i+a_k) \]
The invertible elements of $R_{ij}$ are those which are also in
$R_{ji}$, i.e. $(r^{-1})_{ij}=r_{ji}$.
It follows that  every s lattice homomorphism induces
 a homomorphism of coordinate rings. If $L$ is complemented,
then surjectivity is preserved.

\abs
For a right module $M_R$ let ${\cal L}(M_R)$ denote
 the lattice of all right $R$-submodules.
A {\em von Neumann regular ring} is an associative ring with unit
such that for each $r \in R$ there
is a {\em quasi-inverse}   $x \in R$ such that $rxr=r$
(so homomorphic images are also regular).
Equivalently, the principal right ideals form a complemented  sublattice
 $\ov{L}(R_R)$ of the lattice
${\cal L}(R_R)$ of all
 right ideals - consisting precisely of the
compact elements. And, equivalently, each principal right ideal has an idempotent generator (resp. the same on the left). 
The lattice structure is given  in terms of idempotents $e,f,g$ by
\[ \begin{array}{lclll} eR+fR&=&(e+g)R &\mbox{ with }   & gR=(f-ef)R\\
eR \cap fR&=&(f-fg)R&\mbox{ with }& Rg=R(f-ef)\\
Re \oplus R(1-e)&=&R && \end{array}\]
\begin{corollary} If $R$ is regular and 
 $\phi:R \rightarrow S$ a surjective homomorphism then there is 
  a surjective homomorphism  $\ov{\phi}:\ov{L}(R)
\rightarrow \ov{L}(S)$ such that   $\ov{\phi}(aR)=\phi(a)S$. 
\lab{homo}
\end{corollary}
This is part of  the following result of
Wehrung \cite{Wehr1}

\begin{theorem}  For a regular ring $R$
there is a 1-1-correspondence between two-sided ideals of $R$ and
 neutral ideals of \lab{wehr}
$\ov{L}(R)$ given by
\[ I=\{a\in R\mid aR \in {\cal I}\},\quad\quad {\cal I}=\{aR\mid a \in I\}\]
\end{theorem}
We say that a lattice $L$ is {\em coordinatized} by the regular ring $R$,
if $L$ is isomorphic to $\ov{L}(R_R)$ - and then
 Arguesian, in particular.
Of course, a  height $2$-lattice is coordinatizable
if an only if it is infinite or has $p^k+1$ atoms for some
$k$ and some prime $p$. 
From J\'{o}nsson \cite{J60} Cor.8.5, Lemma 8.2, and Thm.8.3 and von 
Neumann 
\cite{NEU} (see \cite{Hfr} for a short proof)  we have
\begin{theorem}
 Every  complemented modular lattice
 which is simple of height $\geq n$ or has a spanning
frame of order $n$, $n \geq 4$ resp. $n \geq 3$ and $L$ Arguesian,
can be coordinatized by a regular ring.
 Every interval $[0,u]$ of a coordinatizable
lattice is coordinatizable. \lab{J}
\end{theorem}

\abs
\no
We need more information about frames and an alternative view 
of coordinatization. Recall, that the ring $R_n$ of
$n\times n$-matrices over a regular ring $R$ is itself regular.
Assume  $n \geq 3$
and  let ${e}_i$ denote the i-th unit vector in the module $R^n$.
\begin{description}
\item[(1)] Given a   ring $R$, 
the  right submodules of $R^n$ form  a modular lattice ${\cal L}(R^n_R)$.
For regular $R$, the finitely generated ones form
a complemented  sublattice $\ov{L}(R^n_R)$.
 Moreover, 
the $E_i={e}_iR$, $ i \leq n$ and $E_{ij}=({e}_i -{e}_j)R$
form a  spanning ({\em canonical}) frame $\f{E}$. For $n\geq 3$, the lattice
$\ov{L}(R^n_R)$ is generated by $\f{E}$ and its coordinate ring.
\item[(2)] For every complemented modular  $L$ with  spanning  $n$-frame $\f{a}$ 
 there is  regular ring $R$ and an isomorphism $\phi$ of $\ov{L}(R^n_R)$  onto $L$ with 
$\phi(\f{E})= \f{a}$.
Moreover, $R(L,\f{a})_{ij}$ is a regular ring with zero $a_i$, unit
$a_{ij}$, $\oplus_{ij}$, $\ominus_{ij}$, and $\otimes_{ij}$  and an
 isomorphic image of 
 $R$ via $ r_{ij} \mapsto \phi(({e}_i -{e}_jr)R)$.  
\item[(3)] 
The lattices ${\cal L}(R_{nR_n})$ and ${\cal L}(R^n_R)$
 are isomorphic with
an ideal $I$ corresponding to a submodule $U$ iff 
the columns in $U$ are exactly the columns of matrices  in $I$. 
The canonical idempotent matrices with all entries $0$ but one diagonal entry $1$ correspond to the canonical basis vectors.
This isomorphism takes 
$\ov{L}(R_{nR_n})$ to $\ov{L}(R^n_R)$. 
\end{description}

\subsection{Coordinatization of ortholattices}

 An {\em involution} $*$ on a ring $R$ is an involutory anti-automorphism
\[(r+s)^*=r^*+s^*,\quad (rs)^*=s^*r^*,\quad r^{**}=r \quad \mbox{ for all }
r,s \in R.\]
An element such that $r^*=r$ is called {\em hermitian}.
A {\em $*$-ring} is an associative ring $R$ with $1$ endowed with an
involution. 
A $*$-ring is $*$-{\em regular}
  if it is von Neumann regular and if
\[ r^*r=0 \mbox{ implies } r =0 \mbox{ for all } r \in R .\]
Equivalently,  each principal right ideal is generated by an hermitian idempotent.
On a $*$-regular ring $R$,
 $x \perp y \gdw x^*y=0$
 defines an anisotropic symmetric relation 
compatible with addition and right scalar multiplication,
 whence an anisotropic orthogonality on ${\cal L}(R_R)$.
In particular
\[ X \mapsto X\p=\{ y \in R\mid  \forall x \in X.\; x \perp y\} \in {\cal L}(R_R)\]
turns
  $\ov{L}(R_R)$ into an MOL - 
again this characterizes $*$-regularity. This MOL  satisfies  the same
orthoimplications as ${\cal L}(R_R)$ and  is said
to be {\em coordinatized} by $R$. 
If $e$ is a hermitian idempotent we also have
$eR^\perp =(1-e)R$ and $eRe$ is $*$-regular if $e$ is, in addition, central.

\begin{corollary}
 If  $R$ is  $*$-regular and $I$ an  ideal
of $R$ then  $I^*=I$ and $R/I$ is $*$-regular, too.
Homomorphic images of coordinatizable MOLs are  coordinatizable. 
\end{corollary}
{\em Proof.}  $I$ is generated by $\{e\mid e^*=e,\; eR \in {\cal I}\}$,
whence closed under the involution.
Thus, $R/I$ is a $*$-ring, naturally,
and $*$-regular since  every principal right ideal
is generated by a hermitean idempotent.
Thus,
 $\ov{L}(R/I)$  with involution $\p$ is an MOL 
and the lattice homomorphism $\ov{\phi}$,  
associated with  the canonical homomorphism of $R$ onto $R/I$ 
according to Cor.\ref{homo},
preserves orthocomplementation, as well.
The second claim follows by Thm.\ref{wehr}. $\pfbox$
From von Neumann \cite{NEU} II, Thms.4.3 -4.5 and \ref{J}
we have
\begin{theorem}  Every MOL coordinatized as a lattice
by a regular ring is coordinatized by a $*$-regular ring -
having the given ring as reduct. In particular, every MOL
$L$
 with spanning
frame of order $n \geq 4$ ($n \geq 3$ for Arguesian $L$) 
can be coordinatized by a $*$-regular ring. 
\lab{N}
\end{theorem}
Now, assume we are given an  MOL $L$ and $n \geq 3$. 
A frame $\f{a}$ in $L$ is {\em orthogonal}, if $a_j \leq a'_k$ for all $j \neq k$ cf \cite{R:89}.
According to Maeda \cite{MAE,MAE2} we can add to  the above 
description
\begin{description}
\item[(1)] Given a $*$-regular ring $R$ and invertible elements
$\al_1, \ldots ,\al_n$ of $R$ such that $\al_i^*=\al_i$
then
 $L=\ov{L}(R^n_R)$ is a MOL with orthogonal frame $\f{E}$ and 
\[X'=\{(y_1, \ldots ,y_n)\,|\, \sum_{i=1}^n y^*_i\al_ix_i =0 \;\;\mbox{ for all }  
(x_1, \ldots ,x_n) \in X\}.\]
\item[(2)] For every MOL $L$ with  spanning orthogonal  frame $\f{a}$ 
 there is  an isomorphism $\phi$ of an  MOL as in (1) (and w.l.o.g.
 $\al_1=1$) onto $L$ with 
$\phi(\f{E})= \f{a}$.
Moreover, $R(L,\f{a})_{12}$ is $*$-regular 
\item[(3)] The matrix ring  $R_n$  of  $*$-regular ring $R$ is
$*$-regular if and only if there are $\al_i$ as in (1).  
Then, the involution is given by
\[ (x_{ij})^* =( \al^{-1}_ix^*_{ji}\al_j) \] 
and  the isomorphism between
$\ov{L}(R_{nR_n})$ and $\ov{L}(R^n_R)$
is an MOL-isomorphism, too.
\end{description}

\begin{lemma} \lab{matr}
Let $S$ be a $*$-regular ring such that $\ov{L}(S_S)$ contains
an orthogonal $n$-frame $\f{a}$. Then choosing hermitian
idempotents $e_i$ generating $a_i$  there is a $*$ regular ring
$R$ with invertible hermitian $1=\al_1, \ldots ,\al_n$ 
such that $S$ is isomorphic to the $*$-ring $R_n$ as above
and the induced MOL-isomorphism maps $\f{a}$ onto the
canonical frame.
\end{lemma}
{\em Proof.} The case  $n=2$ is illustrative enough.
 We may assume that $S=R_n$ as a ring and 
\[ e_1=E_1=\pmatrix{1&0\cr0&0},\;\;e_2=E_2=\pmatrix{0&0\cr0&1}\] 
To define the involution on $R$ consider
\[ A=\pmatrix{r&0\cr0&0},\;\;A^*=\pmatrix{a&c\cr b&d}\]
Form $A\perp E_2$ we get
$A^*E_2=0$ and $c=d=0$. From $A=AE_1$ we get $A^*=E_1A^*$ and $b=0$.
Thus, we get an involution of $R$ such that
\[ \pmatrix{r&0\cr0&0}^* =\pmatrix{r^*&0\cr 0&0}\]
Using orthogonality to $E_2$ resp. $E_1$ we get
\[\pmatrix{0&1\cr0&0}^*=\pmatrix{a&0\cr \be&0},\;\;\pmatrix{0&0\cr \be&0}^*=\pmatrix{0&c\cr0&d},\;\;
\pmatrix{0&1\cr0&0}=\pmatrix{a&0\cr \be&0}^*
=\pmatrix{a^*&0\cr0&0}+\pmatrix{0&c\cr0&d}\]
whence $a^*=0$ and $a=0$. Thus, with a similar argument,
we have $\be$ and $\al$ in $R$ such that
\[ \pmatrix{0&1\cr0&0}^*=\pmatrix{0&0\cr\be&0},\;\;
\pmatrix{0&0\cr1&0}^*=\pmatrix{0&\al\cr0&0},\;\;
\pmatrix{0&1\cr1&0}^*=\pmatrix{0&\al\cr\be&0}
\]
Hence
\[\pmatrix{0&r\cr0&0}^*=(\pmatrix{r&0\cr0&0}\cdot \pmatrix{0&1\cr1&0})^*
=\pmatrix{0&0\cr \be r^*&0}\]\[
\pmatrix{0&0\cr r&0}^*=(\pmatrix{0&1\cr1&0}\cdot \pmatrix{r&0\cr0&0})^*=
\pmatrix{o&r^*\al \cr 0&0}\]\[
\pmatrix{0&0\cr0&r}^*=(\pmatrix{0&1\cr1&0}\cdot \pmatrix{0&r\cr0&0})^*
=\pmatrix{0&0\cr0&\be r^* \al}
\]
\[\pmatrix{1&0\cr0&1}=\pmatrix{1&0\cr0&1}^*= (\pmatrix{0&1\cr1&0}^2)^*
= \pmatrix{0&\al\cr\be&0}^2=\pmatrix{\al\be&0\cr0&\be\al}\]\[
\pmatrix{0&1\cr1&0}=\pmatrix{0&\al \cr \be &0}^*
=\pmatrix{0& \be^* \al\cr \be \al^* &0} \]
whence $\be=\al^{-1}=\be^*$. $\pfbox$

\abs
Given a right $R$-module $V$,  a map
$\Phi:V^2 \imp R$ is  {\it $*$-sesqui-linear}
if $\Phi(x,y)$ is  linear in $y$,
 additive in $x$, and 
$\Phi(rx,y)=r^*\Phi(x,y)$. It is
 {\it $*$-hermitian} if also $\Phi(y,x)=(\Phi(x,y))^*$. 
Defining
\[  U\p =\{x \in V\mid \forall u \in U.\; \Phi(x,u)=0\}, \quad
L_\Phi(V)=\{ U \in L(V_R)\mid U\h=U \}\]
one obtains a  complete lattice with involution $\p$
which is an ortholattice  if an only
if  $\Phi$ is anisotropic: $\Phi(x,x)\neq 0$ for $x \neq 0$.
If also   $\dim V_R < \infty$ them it is an MOL.  Of course, with respect
to an orthogonal basis, one obtains a description by a diagonal 
matrix as in (1) above.
Now, the results of Baer \cite{RB:46} and Birkhoff and von 
Neumann \cite{BVN:36} can be formulated as follows
\begin{theorem} 
Every finite MOL is a direct product of Boolean algebras and 
MO$_n$'s. Every finite dimensional MOL is a direct product
of MOLs of height $\leq 3$ and 
 MOLs arising from finite dimensional vector spaces
with anisotropic $*$-hermitian form resp. matrix $*$-rings
over skew fields.
\end{theorem}

\section{MOLs in projective spaces}

\subsection{Projective spaces}

If a modular lattice, $M$, is  algebraic (i.e. complete with a join-dense 
set of  compact elements)
and atomistic (equivalently: $M$ is complemented
resp. $1_M$ is a join of atoms) we speak of a {\em
geomodular} lattice.
  By 
$M_{fin}$ we denote the neutral ideal of elements of finite height in $M$.
For  geomodular $M$, 
these are the  elements which are joins of finitely many atoms.

\abs
By a projective space we understand a set $P$ of points together with
a distinguished set of 3-element subsets, the {\em collinear triplets},
 such that the following `triangle axiom' holds: If $p,s,q$ and $q,t,r$ are collinear but $p,q,r$ are not  then
there is unique $u$ such that $p,r,u$ and $s,t,u$  are collinear. 
 A {\em subspace} of $P$ is a subset $U$ of $P$ such that if $p,q \in U$ and
 $p,q,r$ collinear
then $r \in U$. 
The subspaces form a geomodular  lattice ${\cal S}(P)$
where meet is intersection
and the join of $X$ and $Y$ consists of all $r$ collinear with 
some $p \in X$ and $q \in Y$.
 Singleton subspaces and points are identified.
$P$ is {\em irreducible} if for any two points there is a third
one collinear with them. If $P$ is irreducible and ${\cal S}(P)$ 
of  height $n \geq 4$ then, by the Coordinatization Theorem of
 Projective Geometry, there is  a vector space $V$ such that
${\cal S}(P)$ is isomorphic to the lattice ${\cal L}(V)$
of linear subspaces of $V$.

Now, let $M$ be any modular lattice and $P=P_M$ be the set of {\em points},
 i.e. atoms, of $M$.
Then $P$ is turned into a projective space where
 $p,q,r$ are collinear if $p+q=p+r=q+r$.
 We will refer to this as the projective space $P_M$ 
of $M$. The subspace lattice ${\cal S}(P)$
is canonically isomorphic to the ideal lattice
of the sublattice of $L$ consisting of all elements which are
joins of finitely many atoms. 

If $M$ is  algebraic and $P=P_M$, then ${\cal S}(P)$  is isomorphic to the interval sublattice $[0,\sum P]$ of $M$ in the following manner:
If  $u \in M$ then $U =  \{p \in P | p \leq u \}$ is
a subspace.   Conversely, if $S$ is a subspace then
$\sum S \in M$, and $S = \{ p \in P | p \leq \sum S \}$.
  It will sometimes be convenient to consider
$u \in M$ as a subspace and, when we do, we will do so 
without changing notation.
A {\em  subgeometry} $Q$ of a projective
geometry $P$ is just a relatively complemented  $0$-sublattice of
 ${\cal S}(P)_{fin}$ with set $Q \subseteq P$ of atoms. 
In other terms, $Q$ is a subset of $P$ with the induced collinearity
and closed under the operation given by the  triangle-axiom:
If $p,q,r,s,t$ are in $Q$ and $p,s,q$ and $q,t,r$ are collinear, but $p,q,r$ are not, then
there is  $u$ in $Q$ such that $p,r,u$ and $s,t,u$ are collinear.

\abs
The disjoint union of projective spaces $P_i$ constitutes
a projective space $P$. Conversely,
on the point set $P$ of a geomodular lattice, perspectivity is transitive 
and $P$ splits into connected {\em irreducible components} $P_i$ which are
projective spaces in their own right. The subspace lattice ${\cal S}(P_i)$
forms an interval $[0,\sum P_i]$ in ${\cal S}(P)$  and ${\cal S}(P)$ is isomorphic
to the direct product of the ${\cal S}(P_i)$ via
\[ X \mapsto (X\cap P_i\mid i \in I)\]
In particular, we have projections which are lattice homomorphisms
\[ \pi_i:{\cal S}(P) \imp {\cal S}(P_i), \quad \pi_i X = X \cap P_i \]
 The following are due
to  Frink \cite{FRNK:46} (cf \cite{CD}).

\begin{lemma} \lab{frink0}
Let $ab=0$ in $ M$ and $p$ an atom of $M$ such that $p\leq a+b$, $p\not\leq a$, 
and $p\not\leq b$. Then $p, a(p+b), b(p+a)$ are collinear atoms of $M$.
In a complemented modular lattice, if 
$p$ an is atom of $M$ such that $p\leq a+b$, $p\not\leq a$, 
and $p\not\leq b$ then there are atoms $q \leq a$ and $r \leq b$
such that $p,q,r$ are collinear.
\end{lemma}
{\em Proof.} The first is done by a direct calculation. In the
second let $\ti{b}$ a complement of $ab$ in $[0,b]$ and apply
the first. $\pfbox$

\abs
For any map  $\ga:L \imp M$ of  and $Q \subseteq P_M$ there is a natural
map $\ga_Q:L \imp {\cal S}(Q)$ given by
\[ \ga_Q a= \sum\{q \in Q \mid q \leq \ga
a\}\]

\xc{
\setlength{\unitlength}{0.00500000in}
\begingroup\makeatletter\ifx\SetFigFont\undefined%
\gdef\SetFigFont#1#2#3#4#5{%
  \reset@font\fontsize{#1}{#2pt}%
  \fontfamily{#3}\fontseries{#4}\fontshape{#5}%
  \selectfont}%
\fi\endgroup%
{\renewcommand{\dashlinestretch}{30}
\begin{picture}(415,383)(0,-10)
\put(78,218){\blacken\ellipse{10}{10}}
\put(78,218){\ellipse{10}{10}}
\put(198,338){\blacken\ellipse{10}{10}}
\put(198,338){\ellipse{10}{10}}
\put(318,218){\blacken\ellipse{10}{10}}
\put(318,218){\ellipse{10}{10}}
\put(358,178){\blacken\ellipse{10}{10}}
\put(358,178){\ellipse{10}{10}}
\put(238,58){\blacken\ellipse{10}{10}}
\put(238,58){\ellipse{10}{10}}
\put(198,58){\blacken\ellipse{10}{10}}
\put(198,58){\ellipse{10}{10}}
\put(158,58){\blacken\ellipse{10}{10}}
\put(158,58){\ellipse{10}{10}}
\put(198,98){\blacken\ellipse{10}{10}}
\put(198,98){\ellipse{10}{10}}
\put(198,18){\blacken\ellipse{10}{10}}
\put(198,18){\ellipse{10}{10}}
\thicklines
\path(38,178)(198,18)
\path(38,178)(198,18)
\path(198,18)(358,178)
\path(198,18)(358,178)
\thinlines
\put(38,178){\blacken\ellipse{10}{10}}
\put(38,178){\ellipse{10}{10}}
\thicklines
\path(38,178)(198,338)
\path(38,178)(198,338)
\put(0,168){\makebox(0,0)[lb]{\smash{{{\SetFigFont{12}{14.4}{\familydefault}{\mddefault}{\updefault}$a$}}}}}
\path(198,338)(358,178)
\path(198,338)(358,178)
\path(158,58)(318,218)
\path(78,218)(238,58)
\path(198,98)(198,18)
\put(214,0){\makebox(0,0)[lb]{\smash{{{\SetFigFont{12}{14.4}{\familydefault}{\mddefault}{\updefault}$0$}}}}}
\put(203,48){\makebox(0,0)[lb]{\smash{{{\SetFigFont{12}{14.4}{\familydefault}{\mddefault}{\updefault}$p$}}}}}
\put(251,42){\makebox(0,0)[lb]{\smash{{{\SetFigFont{12}{14.4}{\familydefault}{\mddefault}{\updefault}$(p+a)b$}}}}}
\put(372,170){\makebox(0,0)[lb]{\smash{{{\SetFigFont{12}{14.4}{\familydefault}{\mddefault}{\updefault}$b$}}}}}
\put(220,333){\makebox(0,0)[lb]{\smash{{{\SetFigFont{12}{14.4}{\familydefault}{\mddefault}{\updefault}$a+b$}}}}}
\put(38,38){\makebox(0,0)[lb]{\smash{{{\SetFigFont{12}{14.4}{\familydefault}{\mddefault}{\updefault}$(p+b)a$}}}}}
\end{picture}
}
}

\begin{lemma} \lab{frink}
Let $M,L$ be modular lattices, $Q$ a subgeometry of $P_M$,  $L$ complemented, and $\ga:L \imp M$ a 
 lattice homomorphism.
 Then $\ga_Q:L \imp {\cal S}(Q)$ is a lattice homomorphism 
provided  that $\ga 0=0_M$ and
\begin{quote}
 for all $a,b \in L$ with $ab=0$  and all $p \in Q$ with
$p \leq \ga a +\ga b$ but $p\not\leq \ga a$ and $p \not\leq \ga b$
one has  also $(p+\ga b)\ga a \in Q$ and $(p+\ga a)\ga b \in Q$
\end{quote}
Moreover, $\ga_Q$ is a lattice embedding
if $\ga$ is such and
 for
all $a>0$ in $L$ there is $p \in Q$ with $p \leq \ga a$
\end{lemma}

\abs
The {\em Frink embedding} of a complemented modular lattice
 arises by Lemma \ref{frink} from the principal embedding
$\ga_L:L \imp {\cal F}(L)$.
The points are the maximal filters of $L$ - we also 
speak of the {\em Frink space} of $L$.

\abs
If $Q \subseteq P$  is closed under
perspectivity, i.e. a union of components then 
\[ \pi_Q:{\cal S}(P) \imp {\cal S}(Q),\quad \pi_Q x =x \cdot \sum Q=
 \{q \in Q\mid q \leq x\}\]
is a surjective lattice homomorphism.
For any $0$-lattice homomorphism $\vep:L \imp M$  and $Q \subseteq P_M$ there is a natural
map $\vep_Q:L \imp {\cal S}(Q)$ given by
\[ \vep_Q a= \sum\{q \in Q \mid q \leq \vep a\}\]
and this the $0$-lattice homomorphism $\pi_Q \circ \vep$, if
$M={\cal S}(P)$ and  $Q$ is
 closed under perspectivity.
  
Now, let $\vep$ be an embedding - so consider $L$ as a sublattice of
 ${\cal S}(P)$. If $P$ is the disjoint union of 
subspaces $P_i$ then the projections $\pi_i$ provide a subdirect
decomposition of $L$.
Thus, if $L$ is subdirectly irreducible  
then there exists a component $Q$ of $P$ such that $\vep_Q$ is an embedding.

\subsection{Orthogonalities}

By an {\em orthogonality on a lattice} $M$ we understand a binary
relation such that $0 \perp u$ for all $u$ and
\[ u \perp v \mbox{ implies } v \perp u \]
\[ u \perp v \mbox{ and } w \leq v \mbox{ together imply } u  \perp w \]
 \[ u\perp v \mbox{ and } u \perp w   \mbox{ together imply } u \perp v+w\]
The orthogonality {\em induced} on a subset $Q$ of $M$ is given by
\[\mbox{
$x \perp_Q y$ iff $x \perp y$, \quad $x,y \in Q$}\]
Given $M_i$ with orthogonality $\perp_i$, the product $M$ has the orthogonality
\[ (a_i\mid i \in I) \perp (b_i\mid i \in I)\; \mbox{ iff }\; \forall i \in I.\; a_i 
\perp_i b_i \]
If $\phi:M \imp N$ is a surjective homomorphism, then $N$ has 
the orthogonality
\[ a \perp_N b \mbox{ iff }  a=\phi c,\;b=\phi d \mbox{ for some }
 c \perp_M d\]
On the filter lattice ${\cal F}(M)$ we obtain the {\em canonical} orthogonality
\[ F \perp_{\cal F} G \mbox{ iff } a \perp b \mbox{ for some } a \in F,\;b \in G\]
 An orthogonality is {\em anisotropic} if $u \perp v$ implies $uv=0$.
This property is preserved under forming products, sublattices,  homomorphic images, and filter lattices.
An orthogonality is {\em non-degenerate} if $u \perp v$ for all $v$ implies $u=0$.
This is obviously so in the anisotropic case.
  For any ortholattice we have a canonical orthogonality: 
 $x \perp y$,
iff $x \leq y'$.

\abs
Now, let $M$ be algebraic and
 $P$  a join-dense set of compact elements 
   such that for any  $u,v \in M$ 
and $p \in P$ with $p \leq u+v$ there are $q \leq u$ and $r \leq v$ in $P$
such that $p \leq q+r$. This applies with $P$ the
set of all compact elements of any $M$ resp. $P$ the set of points in a geomodular $M$. 
By an {\em orthogonality} on  $P$   we understand a binary  symmetric
relation $\perp$ on $P$ such that
\[ p\perp q,\; p \perp r, \mbox{ and } s \leq q+r \mbox{ together imply } p \perp s\]
We obtain an orthogonality on $M$ defining
 \[ u \perp v \;\;\; \mbox{ iff }  \; p \perp q \mbox{ for all } p \leq u,\, q \leq v \]
Namely, the first two properties are obvious. In the third we may assume $u \in P$. Now, if $v+w \geq p \in P$
then there are
$q \leq v, r \leq w$ such that $p\leq q+r$ whence
$u \perp q$, $u \perp r$ and so $u \perp p$.  
Defining
 \[ u\p =
 \sum\{ q \in P\, |\, q \perp u   \} \]
we get
\[u \leq v\p \mbox{ iff } u \perp v \mbox{ iff }  v \leq u\p \]
Namely, if $u \geq p\in P$ and $u \perp v$, then $v\p \geq p$ and there are finitely many
$q_i \in P$ with $q_i \perp v$ and $p \leq \sum q_i $ whence $p \perp v$.
It follows 
that $x \mapsto x\p$ is a self adjoint Galois connection on the lattice $M$.
  In particular
 the map is order reversing and the map
$x \mapsto x\h$ is a closure operator on $M$. To wit
\[\mbox{$u \leq v$ implies $v\p \leq u\p$,
\quad $u \leq u\h$}\]\[\mbox{
$(\sum_{i\in I}u_i)\p=\prod(u_i\p)$,\quad
  $u\p =
\prod \{p\p\,|\, u \geq p \in P \}$}.\]
The  closed elements of $M$,
endowed with the partial order inherited from $M$ and the restriction of $\p$,
form a
 complete meet sublattice $K$ of $M$ and a complete ortholattice 
containing $L$ as a subalgebra.
Indeed,
$a \vee_K b = a\h \vee_K b\h = (a\p b\p)\p
=(a'b')\p =(a'b')' =a + b$ for $a,b \in L$.
Moreover $K$ is atomistic if $P$ consists of  atoms
$p$ such that $p=p\h$.
$K$ satisfies the {\em covering property}. if
 $u\vee p=u+p$ covers $u$ for any atom $p\not\leq u$.

Let $Q$ be join-dense in $M$.
Then any orthogonality on $M$ is determined by
the orthogonality 
 induced on $Q$.
 It is  anisotropic if $p \not\perp p$ for all $p \in Q$
and non-degenerate if  for
each $p \in Q$ there is $q \in Q$ such that $p \not\perp q$.
For a direct product $M={\cal S}(P)$ of  $M_i={\cal S}(P_i)$ with the product orthogonality,
  we have $P$ the disjoint
union of the $P_i$   and speak of the {\em orthogonal disjoint union} of the $P_i,\perp_i$.

\begin{proposition}
Let $L$ be a bounded lattice with anisotropic orthogonality $\perp$
and assume that for each $x$ there is $x'$ such that
\[ x' \perp x\; \mbox{ and } \forall y.\;\; y \leq x+(x+y)x' \]
Then $x'=\sup\{z\mid z \perp x\}$ is
uniquely determined  and $L$ with $x \mapsto x'$
is an orthomodular lattice. \lab{omo}
\end{proposition}
{\em Proof.} With $y=1$ we get $1 \leq x+x'$. Now, if $x \leq y$ and $y \perp x'$ 
then $y \leq x+yx' =x$ whence $x =\sup\{y \mid y \perp x'\}$.
Since $x''\perp x'$, it follows $x'' \leq x$. For $z \geq x'$ and $z \perp x$
we get $z \leq x'+ x''z\leq x'+ xz = x'$ whence $x'=\sup\{z \mid z \perp x\}$.
It follows $x \leq x''$, thus $x=x''$. Moreover,
$x \leq y'$ implies $x\perp y$ whence $y \leq x'$ and 
we have an ortholattice, indeed. $\pfbox$

\begin{lemma} Every  $0$-$1$-lattice embedding $\eta:L \imp {\cal S}(P)$
of an MOL 
induces an anisotropic  \lab{orthrep} orthogonality on $P$
\[ p \perp q \;\;\mbox{ iff } p \leq \eta a \mbox{ and } q \leq \eta(a')
\mbox{ for some } a \in L\]
Moreover,
\[
\eta (a') \leq  (\eta a)\p \;\; \mbox{ for all } a \in L\]
\end{lemma}
{\em Proof.} For convenience, we think of $\eta$ as $id_L$.
Consider $p \perp q,r$ and $s \leq q+r$.
Then $p \leq a,b$ and $q \leq a, r \leq b$ for some $a,b \in L$
whence $p \leq ab$ and $s \leq a'+b'=(ab)'$. Thus 
we obtain  an anisotropic orthogonality.
Moreover, if $p \leq \eta(a')$ then $p \leq a'$, i.e. $p \perp a$ and
so $p \perp \eta a$. Thus, $\eta(a') \leq (\eta a)\p$.
Now, $\eta a +\eta a'=\eta(a+a')=1_M$
 by embedding, $\eta a \cdot (\eta a)\p=0$ by anisotropicity,
and $\eta a' \leq  (\eta a)\p$ by hypothesis, whence $\eta a' =(\eta a)\p$
by modularity. $\pfbox$

\subsection{Polarities}

An orthogonality on a geomodular lattice resp. its projective space is a {\em polarity} if it is nondegenerate and if $p\p$ is a coatom  for each atom $p$.

\begin{lemma} A nondegenerate  orthogonality $\perp$ on a geomodular lattice
is a polarity if and only if
\[ p\p(q+r)>0 \;\;\mbox{ for all points } q\neq r 
\mbox{ with } p \not\perp q,\;p \not\perp r\] 
\end{lemma}
{\em Proof.}  If $\perp$ is a polarity then  $p\p$ is a coatom whence 
the claim follows by modularity.
Conversely, consider $q \not\leq p\p$. We have to show that $q+p\p=1$,
i.e. $r \leq q+p\p$ for all $r \neq q$. But, by hypothesis,
if $r \not\leq p\p$ then $0<s=p\p(q+r)<q+r$, so $s \leq p\p$ is a point
and $r \leq q+s \leq q+p\p$. $\pfbox$

\begin{corollary} For an anisotropic orthogonality $\perp$ on a geomodular lattice
the following are equivalent
 \bd
\item[(1)] $\perp$  is a polarity 
\item[(2)]  $p+p\p=1$ for all atoms $p$
\item[(3)] $(p+r)p\p>0$ for all atoms $p\neq r$ \label{pol1}\note{pol1}
\ed\end{corollary}
{\em Proof.}  As observed above, $p p\p=0$. Thus for $r \neq p$ we have
(3) trivially, if $r \perp p$, and by the Lemma, otherwise.
If (2) holds, then $p\p$ is a coatom by modularity.
 Thus, by modularity, $p+p\p=1$ if and only 
if $p\p$ is a coatom. $\pfbox$ 
\begin{corollary} For each MOL $M$ there is a canonical
anisotropic polarity on $P_M$ given by
\[ p \perp q \;\mbox{ if and only if } \; p \leq p'\]
\end{corollary}

\begin{corollary}
A projective space with anisotropic polarity is the orthogonal
disjoint union  of its irreducible
components. Conversely, the orthogonal disjoint union of spaces
with polarity yields a space with polarity.\lab{odec}
\end{corollary}
{\em Proof.}
In view of (3) $p\neq r$ and $p \not\perp r$
jointly imply that there is a  $q$ collinear with $p,r$. $\pfbox$
According to Maeda \cite{MAE} an orthogonality $\perp$ on
a desarguean 
irreducible projective space $P$
(so $L(P) \cong L(V_{D})$ for some vector space)
is a  polarity  if and only if  there is an
anti-automorphism $*$ of $D$ and    $*$-hermitian form
$\Phi$ on $V_D$  such that
\[ p=vD \perp q=wD \;\mbox{ if and only if }  
\Phi(v,w)=0 \]
and  $\Phi$ is anisotropic if and only if so is $\perp$. 
For such, the lattice $L_\Phi(V)$ of closed elements is modular
if and only if $V_D$ is finite dimensional
(Keller \cite{KEL}).

\begin{lemma} Let $\perp$ be an anisotropic polarity on the geomodular
lattice $M$. Then 
 \bd
\item[(1)] $u\h=u$ and $u+u\p=1$ for all $u \in M_{fin}$
\item[(2)] Each interval $[0,u] \subseteq M_{fin}$ with the induced
orthogonality  is an MOL with orthocomplementation
\[ x \mapsto x^{\perp u} = \sum \{q \leq u\,| \, q \perp x \}=
u x\p.\]
  \label{fin}\note{fin}
\ed
\end{lemma} 
{\em Proof} by induction on the height of $u$. For $u=0$ nothing is to be done.
So let $v$ a lower cover of $u$. By inductive hypothesis.
$v+v\p=1$, whence by modularity $p=uv\p \in P$ and $u=v+p$. It follows, with modularity again, $u+u\p=v+p +v\p p\p = (v+v\p)(p+p\p)=1$. Since $u\h \geq u$
and $u\p u\h=$ we have $u\h=u$ by another application of modularity.
Finally,
choose $p \leq x$ and let $v=up\p$. $v$ is a lower cover of $u$. Then 
\xc{
\parbox{8cm}{\hspace*{1cm}\setlength{\unitlength}{0.00500000in}
\begingroup\makeatletter\ifx\SetFigFont\undefined%
\gdef\SetFigFont#1#2#3#4#5{%
  \reset@font\fontsize{#1}{#2pt}%
  \fontfamily{#3}\fontseries{#4}\fontshape{#5}%
  \selectfont}%
\fi\endgroup%
{
\begin{picture}(599,451)(0,-10)
\put(460,298){\blacken\ellipse{10}{10}}
\put(460,298){\ellipse{10}{10}}
\put(420,258){\blacken\ellipse{10}{10}}
\put(420,258){\ellipse{10}{10}}
\put(500,138){\blacken\ellipse{10}{10}}
\put(500,138){\ellipse{10}{10}}
\put(100,138){\blacken\ellipse{10}{10}}
\put(100,138){\ellipse{10}{10}}
\put(20,258){\blacken\ellipse{10}{10}}
\put(20,258){\ellipse{10}{10}}
\put(60,298){\blacken\ellipse{10}{10}}
\put(60,298){\ellipse{10}{10}}
\put(260,178){\blacken\ellipse{10}{10}}
\put(260,178){\ellipse{10}{10}}
\put(220,138){\blacken\ellipse{10}{10}}
\put(220,138){\ellipse{10}{10}}
\put(300,18){\blacken\ellipse{10}{10}}
\put(300,18){\ellipse{10}{10}}
\put(340,58){\blacken\ellipse{10}{10}}
\put(340,58){\ellipse{10}{10}}
\put(220,378){\blacken\ellipse{10}{10}}
\put(220,378){\ellipse{10}{10}}
\put(260,418){\blacken\ellipse{10}{10}}
\put(260,418){\ellipse{10}{10}}
\put(340,298){\blacken\ellipse{10}{10}}
\put(340,298){\ellipse{10}{10}}
\put(300,258){\blacken\ellipse{10}{10}}
\put(300,258){\ellipse{10}{10}}
\put(140,178){\blacken\ellipse{10}{10}}
\put(140,178){\ellipse{10}{10}}
\thicklines
\path(420,258)(500,138)
\path(500,138)(540,178)
\path(420,258)(460,298)
\path(460,298)(540,178)
\path(20,258)(100,138)
\path(100,138)(140,178)
\path(20,258)(60,298)
\path(60,298)(140,178)
\path(220,138)(300,18)
\path(300,18)(340,58)
\path(220,138)(260,178)
\path(260,178)(340,58)
\path(220,378)(300,258)
\thinlines
\put(540,178){\blacken\ellipse{10}{10}}
\put(540,178){\ellipse{10}{10}}
\thicklines
\path(300,258)(340,298)
\put(508,121){\makebox(0,0)[lb]{\smash{{{\SetFigFont{8}{9.6}{\familydefault}{\mddefault}{\updefault}$u\p$}}}}}
\path(220,378)(260,418)
\path(260,418)(340,298)
\path(300,18)(500,138)
\path(340,58)(540,178)
\path(260,178)(460,298)
\path(220,138)(420,258)
\path(140,178)(340,298)
\path(100,138)(300,258)
\path(60,298)(260,418)
\path(20,258)(220,378)
\path(100,138)(300,18)
\path(140,178)(340,58)
\path(300,258)(500,138)
\path(340,298)(540,178)
\path(260,418)(460,298)
\path(60,298)(260,178)
\path(20,258)(220,138)
\put(30,298){\makebox(0,0)[lb]{\smash{{{\SetFigFont{8}{9.6}{\familydefault}{\mddefault}{\updefault}$u$}}}}}
\put(0,247){\makebox(0,0)[lb]{\smash{{{\SetFigFont{8}{9.6}{\familydefault}{\mddefault}{\updefault}$v$}}}}}
\put(341,74){\makebox(0,0)[lb]{\smash{{{\SetFigFont{8}{9.6}{\familydefault}{\mddefault}{\updefault}$p$}}}}}
\put(557,174){\makebox(0,0)[lb]{\smash{{{\SetFigFont{8}{9.6}{\familydefault}{\mddefault}{\updefault}$v\p$}}}}}
\put(477,290){\makebox(0,0)[lb]{\smash{{{\SetFigFont{8}{9.6}{\familydefault}{\mddefault}{\updefault}$(vx)\p$}}}}}
\put(318,0){\makebox(0,0)[lb]{\smash{{{\SetFigFont{8}{9.6}{\familydefault}{\mddefault}{\updefault}$0$}}}}}
\put(286,414){\makebox(0,0)[lb]{\smash{{{\SetFigFont{8}{9.6}{\familydefault}{\mddefault}{\updefault}$1$}}}}}
\put(435,253){\makebox(0,0)[lb]{\smash{{{\SetFigFont{8}{9.6}{\familydefault}{\mddefault}{\updefault}$x\p$}}}}}
\put(236,371){\makebox(0,0)[lb]{\smash{{{\SetFigFont{8}{9.6}{\familydefault}{\mddefault}{\updefault}$p\p$}}}}}
\put(106,174){\makebox(0,0)[lb]{\smash{{{\SetFigFont{8}{9.6}{\familydefault}{\mddefault}{\updefault}$x$}}}}}
\end{picture}
}
}}
\[ x=p+vx,\;\;x\p=p\p(vx)\p,\;\; x^{\perp u} =(vx)^{\perp v} \]
whence by induction
\[ x+ x^{\perp u} =p+vx +x^{\perp v} =p+v =u . \quad \pfbox\]
A {\em geometric representation} of an MOL is a $0$-$1$-lattice embedding
$\eta:L \imp M={\cal S}(P)$ into the subspace lattice of
a projective space $P$ with anisotropic polarity $\perp$ such that
\[ \eta(a) \perp \eta(a') \quad \mbox{ for all } a \in L \]
By modularity it follows
\[ \eta(a') =\eta(a)\p \quad \mbox{ for all } a \in L\]
Indeed, $\eta(a)\p \geq \eta(a')$, $\eta(a) \cdot \eta(a)\p=0$,
and $\eta(a)+\eta(a')=\eta(a+a')=1$.

\begin{corollary} \lab{canrep}
Every subalgebra $L$ of an  atomic MOL $M$  has a geometric representation $\eta:L \rightarrow {\cal S}(P_M)$
with $\eta(a)=\{p \in P_M\mid p \leq a\}$.
\end{corollary}

\subsection{Geometric MOL construction}

For each polarity on a geomodular lattice $M$ the following hold
\begin{itemize}
\item[(i)] If $x \leq y \in M$ such that $\dim y/x<\aleph_0$ then
$\dim x\p/y\p \leq \dim y/x$
\item[(ii)] If $x \leq y \in C$ such that $\dim y/x<\aleph_0$ then
 $\dim_C y/x=\dim x\p/y\p =\dim y/x$
\item[(iii)] If $u \in C$ and $x \geq u$ in $M$ such that 
$\dim x/u < \aleph_0$ then $x \in C$ \end{itemize}
Namely, consider $x \leq y$ in $M$ with  $\dim y/x< \aleph_0$. Then  $y=x+\sum p_i$ with $\dim y/x$ many $p_i\in P$ and
$y\p=x\p  \prod_i p\p_i$. This proves (i).
Now, if $x,y \in C$ then 
$\dim y/x =\dim x\p/y\p$.
If $x\prec_C y$ is a covering in $C$, then $y\p < x\p$ and we may choose $p \leq x\p$,   $p \not\leq y\p$. Then $y \not\leq p\p$ and
$y p\p \in C$. It follows
$x=y p\p \prec_M y$ whence (ii).
Finally, if $u \in C$ and $x \geq u$ in $M$ then 
$ \dim x/u \geq \dim u\p/x\p \geq \dim x\h/u\h
=\dim x\h/u\geq \dim x/u$. Thus (iii).

\abs
An important  congruence relation $\mu$ on any modular lattice $M$ (cf \cite{CD})
is given by 
\[\begin{array}{rcl} x\, \mu \, y &\mbox{ iff } &
\dim (x+y)/(x  y) < \aleph_0\\
&\mbox{ iff } & \dim z/x<\aleph_0 \mbox{ and } \dim y/z<\aleph_0
\mbox{ for some } z \geq  x,y\\&\mbox{ iff } & \dim x/u<\aleph_0 \mbox{ and } \dim y/u<\aleph_0
\mbox{ for some } u \leq x,y \end{array}\]
Given any subset $L$ of $M$ we define
\[\hat{L}= \{ x \in C\mid x \,\mu\,u \mbox{ for some } u \in L\}\]
Consider the conditions
 \[\begin{array}{lrl}
(a)& a b\;  \in \hat{L} & \mbox{ for all } a,b \in L\\
(b)& a+b\;  \in \hat{L},\;a\p +b\p \in C& \mbox{ for all } a,b \in L\\
(c)& a\p \in \hat{L} &\mbox{ for all } a \in L
\end{array} \]

\begin{lemma}
Let $\perp$ be a polarity on the geomodular
lattice $M$. 
Then 
\begin{itemize}
\item $(a)$ implies that $\hat{L}$ is  meet-closed in $M$ and $C$, 
simultaneously
\item $(b)$  implies that $\hat{L}$ is  join-closed in $M$ and $C$, simultaneously 
\item $(c)$ implies that $\hat{L}$ is closed under $x \mapsto x\p$ 
 \end{itemize}
In particular, $\hat{L}$ is a modular ortholattice if
$\perp$ is anisotropic and $(a),(b),(c)$ hold. \lab{hat}
\end{lemma}
This is basically Lemma 2 of \cite{BR}. {\em Proof}.
Observe that 
\[\begin{array}{rcl} \hat{L}&=&\{x \in C\,|\, \exists a \in L.\, \exists y, z \in C. \; y \leq z,\;
 a,x\in[y,z] \mbox{ and } \dim z/y<\aleph_0 \}\end{array}\]
In particular,
\[ x,y,z \in C,\; a \in \hat{L},\;  y\leq z,\;
 a,x\in[y,z], \mbox{ and } \dim z/y<\aleph_0 \quad \mbox{jointly imply }\; x \in \hat{L}\] 
Indeed, for $x\,\mu\, a$ in $C$  we have also $y=x a$ and, by (iii), $z=x+a$ in 
$C$ and $y\,\mu\,z$.
Assuming $(c)$, for $x \in \hat{L}$ with (i)  we conclude $y\p \,\mu\,z\p$ whence
$x\p\, \mu\, a\p$ and so $x\p \in \hat{L}$.
 
Now, 
 consider $y\leq a \leq z$ and $v \leq b \leq  w$ in $C$,
$\dim z/y<\aleph_0$, and $\dim w/v<\aleph_0$. Let 
 $x \in [y,z]$ and $u \in [v,w]$. 
By the congruence properties of $\mu$ one has     $x  u \,\mu \,a  b$
and $x+u \,\mu\,a+b$.
By  $(iii)$ $x,u,x u \in C$. Thus $x  u \in \hat{L}$
if $a,b \in L$ and (a).
Moreover $x+u \in C$
provided that $x\geq a,\,u\geq b$ and $a+b \in C$. 

Now  suppose (b) and $a,b \in L$. We show $y+v \in C$ by induction on
$\dim a/y+\dim b/v$. In doing so, by (iii) we may assume that we have
$y \prec t \leq a$ with $t$ and $t+v$ in $C$. Considering the sublattice
of $M$ generated by $y,t,v$ two cases are possible: 
firstly, $y+v=t+v$ with nothing left to do; secondly, 
$y+v \prec t+v$. 
If we had $v\p  y\p \leq t\p$ then by modularity
$v\p +t\p < v\p+y\p$. Now $a\p \leq t\p \leq y\p$, $b\p \leq v\p$
and $a\p+b\p \in C$ by hypothesis.
Thus, as shown above,  we would have $v\p+t\p$ and $v\p +y\p$ in $C$.
It would follow $v t=(v\p+t\p)\p<(v\p+y\p)\p=v y$, a contradiction. So we may
choose $p \in P$ such that $p \leq v\p y\p$, $p \not\leq t\p$. 
Then $p\p \geq y+v$, $p\p \not\geq t+v$. Consequently,
$y+v=(t+v) p\p \in C$.
 With (iii)   it follows
$x+u \in C$ for all $x \in [y,z]$, $u \in [v,w]$
whence $x+u \in \hat{L}$ since $a+b \in \hat{L}$ by hypothesis. 
 $\pfbox$

\begin{proposition} For any geometric representation $L \subseteq M$
of an MOL, 
there is a  sub-MOL $\hat{L}$  of
the ortholattice $K$ of closed elements of $M$ containing $L$ and all atoms of $M$. In particular, $\hat{L}$ is an atomic MOL containing
$L$ as a sub-MOL.
   \lab{geovar}
\end{proposition}
{\em Proof}. Apply \ref{hat}. $\pfbox$
The original example in \cite{BR} was based on a
 separable real Hilbert space $(H,\Phi)$ and 
$L=\{0,H,A,A\p,C,C\p,D,D\p\} \subseteq L_\Phi(H)$
such that  $A\p+C \in L_\Phi(H)$ coatom,  $X+Y=H$ for $X \neq Y$
in $L\setminus \{0\}$, else.
Thus 
$\hat{L}/\,\theta_{fin} \cong MO_3$
whence  $\hat{L}$ is not coordinatizable.
On the other hand, $\hat{L}$ 
contains an  infinite set of orthogonal perspective elements
and is
not  normal  in the sense of Wehrung \cite{Wehr2}.
The same holds for the subalgebra generated by $A,C,D$.
\xc{
\parbox{8cm}{\hspace*{.2cm}\setlength{\unitlength}{0.00500000in}
\begingroup\makeatletter\ifx\SetFigFont\undefined%
\gdef\SetFigFont#1#2#3#4#5{%
  \reset@font\fontsize{#1}{#2pt}%
  \fontfamily{#3}\fontseries{#4}\fontshape{#5}%
  \selectfont}%
\fi\endgroup%
{\renewcommand{\dashlinestretch}{30}
\begin{picture}(598,380)(0,-10)
\put(278,340){\blacken\ellipse{10}{10}}
\put(278,340){\ellipse{10}{10}}
\put(278,260){\blacken\ellipse{10}{10}}
\put(278,260){\ellipse{10}{10}}
\put(158,140){\blacken\ellipse{10}{10}}
\put(158,140){\ellipse{10}{10}}
\put(278,100){\blacken\ellipse{10}{10}}
\put(278,100){\ellipse{10}{10}}
\put(278,20){\blacken\ellipse{10}{10}}
\put(278,20){\ellipse{10}{10}}
\put(398,140){\blacken\ellipse{10}{10}}
\put(398,140){\ellipse{10}{10}}
\put(518,180){\blacken\ellipse{10}{10}}
\put(518,180){\ellipse{10}{10}}
\put(158,220){\blacken\ellipse{10}{10}}
\put(158,220){\ellipse{10}{10}}
\thicklines
\put(158,220){\ellipse{56}{56}}
\put(38,180){\ellipse{56}{56}}
\put(398,220){\ellipse{56}{56}}
\put(398,140){\ellipse{56}{56}}
\put(518,180){\ellipse{56}{56}}
\put(278,60){\ellipse{80}{120}}
\put(278,300){\ellipse{80}{120}}
\put(158,140){\ellipse{56}{56}}
\thinlines
\put(398,220){\blacken\ellipse{10}{10}}
\put(398,220){\ellipse{10}{10}}
\thicklines
\path(278,340)(158,220)
\path(158,220)(278,100)
\thinlines
\put(38,180){\blacken\ellipse{10}{10}}
\put(38,180){\ellipse{10}{10}}
\thicklines
\path(278,340)(398,220)
\put(127,202){\makebox(0,0)[lb]{\smash{{{\SetFigFont{12}{14.4}{\familydefault}{\mddefault}{\updefault}$A$}}}}}
\path(398,220)(278,100)
\path(278,260)(158,140)
\path(278,260)(398,140)
\path(158,140)(278,20)
\path(398,140)(278,20)
\path(278,100)(278,20)
\path(278,340)(278,260)
\path(278,340)(38,180)
\path(278,340)(518,180)
\path(518,180)(278,20)
\path(38,180)(278,20)
\put(0,167){\makebox(0,0)[lb]{\smash{{{\SetFigFont{12}{14.4}{\familydefault}{\mddefault}{\updefault}$D$}}}}}
\put(404,202){\makebox(0,0)[lb]{\smash{{{\SetFigFont{12}{14.4}{\familydefault}{\mddefault}{\updefault}$C\p$}}}}}
\put(405,145){\makebox(0,0)[lb]{\smash{{{\SetFigFont{12}{14.4}{\familydefault}{\mddefault}{\updefault}$A\p$}}}}}
\put(527,166){\makebox(0,0)[lb]{\smash{{{\SetFigFont{12}{14.4}{\familydefault}{\mddefault}{\updefault}$D\p$}}}}}
\put(294,330){\makebox(0,0)[lb]{\smash{{{\SetFigFont{12}{14.4}{\familydefault}{\mddefault}{\updefault}$1$}}}}}
\put(288,1){\makebox(0,0)[lb]{\smash{{{\SetFigFont{12}{14.4}{\familydefault}{\mddefault}{\updefault}$0$}}}}}
\put(125,132){\makebox(0,0)[lb]{\smash{{{\SetFigFont{12}{14.4}{\familydefault}{\mddefault}{\updefault}$C$}}}}}
\end{picture}
}
}
}

\subsection{Topological  MOL construction}
In his paper \cite{FRNK:46} Frink pointed out that his embedding
can be seen as a generalization of Stone's representation of Boolean
algebras as rings of sets.  
In \cite{JON:54} J\'{o}nsson established as much of a duality
as appears possible without an orthogonality.
 Topological representations  for orthomodular lattices
have been given by Iturrioz \cite{It1,It2}. But
 modularity hardly can be characterized within  that approach. Therefore,
we prefer to work on a projective space at the price
of using a more general concept of `topology', as explained in 
Abramsky and Jung \cite{JUNG}.

\abs
An abstract characterization of the Frink embedding    has been given by
J\'{o}nsson \cite{JON:54}: Considering  $L$ as a sublattice  of $M={\cal S}(P)$
it is a {\em regular sublattice} 
which means that
$L$ is a complemented $0$-$1$-sublattice of the geomodular
lattice $M$ such that \bd \item for all 
 $X \subseteq L$ with  $0 =\prod_M X$ then there is finite $Y \subseteq X$
 with $\prod_M Y=0$  \item for any
 $u\in M_{fin}$ and $q \in P$  with $uq=0$ there are $a,b \in L$
with  $a \geq u$,  $b \geq q$, and
 $ab=0$. \ed

\dale
\no
A {\em subspace topology} ${\cal O}$  on a projective space $P$
is a $0$-$1$-sublattice of ${\cal S}(P)$ closed under
arbitrary joins. The members of ${\cal O}$ are referred to as {\em open subspaces}. The space is {\em strongly Hausdorff} if for any finite $n$ and
$p \neq q_i, (1 \leq i\leq n$  in $P$ there are
$U,V \in {\cal O}$ such that $p \in U$, $q_i \in V (i \leq n)$ and $U \cap V=\emptyset$. The space is {\em Hausdorff} if this holds for $n=1$.

An {\em s-basis} ${\cal B}$ of ${\cal O}$ is a $0$-sublattice
such that each member of ${\cal O}$ is a directed sum (i.e. union)
of members of ${\cal B}$. 

Call a subspace $A$  {\em s-compact} if for any covering
$A \subseteq \bigcup_{i \in I} U_i$ with a directed system of open subspaces $U_i$ there
is  $i\in I$ such that $A \subseteq U_i$.
Observe that if $U,V$ are s-compact subspaces then so is $U+V$.

A {\em MOL-space} is a  projective space $P$ endowed with an 
anisotropic orthogonality $\perp$ and a
s-compact  subspace topology ${\cal O}$
having
 a s-basis ${\cal B}$ such that
 $U\p \in {\cal O}$ 
and $U+U\p=P$ for all $U \in {\cal B}$.

If the collinearity relation on $P$ is empty, then $U+V =U\cup V$ and
$U\p=P \setminus U$ which means that in this case MOL-spaces are just Boolean      spaces.

\begin{proposition}  A MOL-space $P$ 
has a unique s-basis, namely the s-compact open subspaces. 
These form a subalgebra $L$ of $({\cal S}(P),\p)$ which is an  MOL.
If $\perp$ is a polarity, the Hausdorff
property implies its strong variant.
\end{proposition}
{\em Proof.}
 Let $X$ be a subspace of a MOL-space $P$ such that
$X$ and $X\p$ are open and $X+X\p=P$. Then $X$ is s-compact and
$X=X\h$.
Namely, let $X= \bigcup U_i$ and $X\p =\bigcup V_j$
directed unions of basic sets, each including $\emptyset$.
 Then $P=\bigcup (U_i+V_j)$ is also
a directed union of basic sets. S-compactness of $P$ yields that
$P=U_i+V_j$ for some $i,j$. By $X\h \cap X\p=0$ and modularity, one derives
 $U_i=X=X\h$. 

It follows that the basic sets are s-compact open - the converse being trivial.
Also, if $U$ is basic, then $U=U\h$. Thus,
  applying the above to 
$X=U\p$ and  $X\p=U$ we get that $X\p$ is s-compact whence basic.
In particular, $L={\cal B}$ is an MOL.

Now, assume $\perp$ a polarity.
For $u \in M_{fin}$ and $q\in P$ with $p \perp u$ there is $a \in L$
such that $u \leq a$ and $q \leq a'$. We show this by induction on the
 height of $u$. For $u=0$ this is trivial. So let $u>0$ and $v$ a lower cover
 of $u$. Then $p=uv\p \in P$ and $p \perp v$. Hence, by inductive
hypothesis we have $a \in L$ such that $a \geq v$
and $a\p \geq p$. Since $p \perp q$ we have $b \in {\cal O}$  such 
that $b \geq p$ and $b\p \geq q$. Since $L$ is a basis, we may 
choose $b \in L$. Then $a+b \in L$ with $a+b \geq p+v=u$
and $(a+b)\p =ap\ b\p \geq q$.

Consider  $0<u \in M_{fin}$ and $q \in P$ with $uq=0$.  Then $v=uq\p$ is a lower cover of $u$ whence $p=uv\p \in P$
and $v \perp q$ as well as $v \perp p$.  
As just shown, we have $a,b \in L$ such that $a,b \geq v$, $a\p \geq q$ and $b\p \geq p$ and we may assume $a \leq b$ and $b\p \leq a\p$. By the Hausdorff property we 
have $c,d \in L$ such that $p \leq c$, $q \leq d$ and $cd =0$.
We may assume $c \leq b\p$ and $d \leq a\p$.
It follows $u=p+v \leq a+c$ and, by modularity, $(a+c)d \leq (a+b\p)a\p=b\p$
whence $(a+c)d =(a+c)b\p d= (ab\p+c)d =cd =0$.
$\pfbox$

An {\em MOL-space} is {\em Frinkian}
 if it is strongly Hausdorff and if  $P=U_i$ for some $i \in I$ 
whenever $P =(\bigcup_{i \in I} U_i)\h$ for a directed system
of open subspaces.

\begin{theorem}  Frink spaces of MOLs with canonical orthogonality and
 basic open
subspaces \[U(a)=\{ p \in P\,|\,p \leq \vep a\}, \quad a \in L.\]
are  Frinkian MOL-spaces. Moreover,
$a \mapsto U(a)$  provides an (object)duality between
  MOLs and Frinkian  MOL-spaces.
\end{theorem}
{\em Proof.}
Consider a Frinkian MOL-space,
By the Proposition,  $L={\cal B}$ is a MOL. 
We claim that $L$ is a regular
sublattice of $M={\cal S}(P)$. If we have $a_i \in L$ such that
$\prod_{i \in I} a_i =0$ then $(\sum_{i \in I} a_i\p)\p=0$ 
and $P=(\sum_{i \in I} a_i\p)\h$. Hence 
$P= \sum_{i \in J}a_i\p$ for some finite $J \subseteq I$
and $0= \prod_{i \in J} a_i$.

\dale
Conversely, let $M$ be the Frink-extension of the MOL $L$.
Then the  $U(a)$, $a \in L$ form a s-basis of s-compact open subspaces.
Namely, observe that $U(a)\p =U(a')$ and
 suppose
that a directed set $\{a_i \in L\,|\,i \in I\}$ is given such that $U(a)=\bigcup_{i \in I} U(a_i)$.
Then in $M$ we have
 $a =\sum_{i\in I} a_i$. Also $\prod_{i \in I}aa'_i=a(\sum_{I\in I} a_i)\p=0$.
 Thus, by  regularity there is
 $j\in I$ with $aa'_j=0$. It follows $a'+a_j =1$ 
whence  $a= a_j$ by modularity and $a \leq a_i$. 

Similarly, if we have $P=(\bigcup_{i \in I} U(a_i))\h$ with 
directed $a_i \in L$ then $0 =\prod_{i \in I} a'_i$
whence, by regularity of the embedding,
 $0=a'_i$ for some $i$ and so $P=U(a_i)$.

Regularity implies the strong Hausdorff property, immediately.
Also 
 if $p \perp q$
then $p \in U(a)$ and $q \in U(a')$ for some $a$.

This shows that we have a Frinkian MOL-space, indeed, and
that  $a \mapsto U(a)$ is an isomorphism of $L$
onto the algebra of s-compact open subspaces.

On the other hand, starting with a Frinkian MOL-space $P$,
as we have seen above, the embedding of $L$ into ${\cal S}(P)$
is regular and  Thm. 2.6 of J\'{o}nsson \cite{JON:54}
applies to show that
\[ \psi(p) =\{a \in L\,|\, p \in a \},\quad \psi(x)= \sum \{\psi p\,|\, p \leq x\}\]
is a lattice isomorphism of ${\cal S}(P)$ onto the
subspace lattice of the Frink-space   
such that $\psi|L$ is the Frink-embedding.  Moreover,
in $P$ we have, by hypothesis, $p \perp q$ iff $p \leq a$ and
$q \leq a\p$ for some basic $a$. Thus, $\psi$ is also
an isomorphism  with respect to  orthogonality.
Since it matches bases  
and it  is a homeomorphism, indeed.
$\pfbox$

\abs
\small Let us take the opportunity to point out an error in A.Day
and C.Herrmann, Gluings of modular lattices, Order {\bf 5} (1988), 85-101.
It is claimed there that the direct limits of the lattices
$({\cal IF})^n(L)$  resp. $({\cal FI})^n(L)$ (taken over the
canonical embeddings) are isomorphic - here ${\cal I}(L)$
denotes the ideal lattice.  Yet, the map $\al$ offered,
 fails to be an isomorphism - and we suspect that there is none. Nethertheless their Lemma 2.1
can be proved directly.

\normalsize
\section{Equational theory}

\subsection{Orthoimplications and varieties}

Let $M$ be a  lattice with $0$ and an orthogonality
(actually, for the generalities we only need  that 
$a \perp b$, $c \leq a$, and $d \leq b$ imply $c \perp d$). 
Considering $M$ as  structure
$(M;+,\cdot,0,\perp)$, the {\em orthoimplication}
given by a lattice term $f$ (in two sorts of variables, $x_i$ and $y_i$)  is the first order
formula,
\[  x_{1} \perp y_{1} \wedge ... \wedge x_{n} \perp y_{n}
\rightarrow f(x_{1},y_{1},...,x_{n},y_{n}) = 0.\]

\begin{lemma} Orthoimplications are preserved under
formation of direct unions,  products, sublattices, homomorphic images, 
and filter lattices - with the induced orthogonalities.
Also, they are preserved under weakening of the orthogonality. 
\lab{HSP}\end{lemma}
{\em Proof}. Formation of direct unions, products, and  substructures 
(weak with respect to the relation symbols) preserves any universal sentences of the above type. Now, let $\phi:L \imp M$ a surjective homomorphism.
Assume $a_i \perp b_i$ in $M$. Then there are $c_i \perp d_i$ in $L$
 with $a_i=\phi c_i$ and $b_i=\phi d_i$. By hypothesis
$f(c_1,d_1, \ldots ,c_n,d_n)=0$ whence 
$f(a_1,b_1, \ldots, a_n,b_n)=0$.
If $F_i \perp_{\cal F} G_i$ then $a_i \perp b_i$ for
some $a_i \in F_i, b_i \in G_i$ whence
\[0=f(a_1,b_1, \ldots ,a_n,b_n) \in  f(F_1,G_1, \ldots ,F_n,G_n) \quad \pfbox\]

\begin{lemma} Let $M$ be an algebraic lattice and $I$ its set of
compact elements or $M$ a complemented modular lattice and
 $I$ a  neutral ideal. 
For each lattice polynomial
$f(z_{1},...,z_{m})$ with constants in $M$ and $c_{1},...,c_{m}$ in $M$,
and for each $p \in I$ one has: 
$f(c_{1},...,c_{m})\geq p$   iff $f(u_{1},...,u_{m}) \geq 
p$ , 
for some  some   $u_{i} \in I$ with $u_{i} \leq  c_{i}$,
$i = 1, ..., m$. 
\lab{22}
\end{lemma}
{\em Proof.}
We proceed by induction on the complexity of
$f$.  The claim is trivially true if
$f$ is a single variable or constant.

Suppose $f = f_{1}  f_{2}$.  Then $p \leq f( c_{1}, ...,  c_{m})$
implies $p \leq f_{k}( c_{1}, ..., c_{m})$, for $k = 1,2$.  By the
inductive hypothesis, there exist  $u_{k1},..., u_{km} \in I$ 
 with $u_{ki} \leq  c_{i}$, $i = 1, ..., m$, and
$p \leq f_{k}(u_{k1},..., u_{km})$, for $k = 1,2$.  Set
 $u_{i} 
= u_{1i} + u_{2i}$, for $i = 1,...,m$.

Now suppose $f = f_{1} + f_{2}$, and, for convenience, define $d_{k}
= f_{k}( c_{1},..., c_{m})$, for $k = 1,2$. 
 In the first case,
we have $d_i =\sum Q_i$ with directed $Q_i \subseteq I$ whence 
by compactness $p \leq \sum P$ with finite $P \subseteq Q_1 \cup Q_2$
and $p_i=\sum P\cap Q_i \in I$.
 In the second case
let $p_i=d_i(p+d_j)$ and $q_i$ a complement of $d_ip$ in 
$[0,p_i]$ and $q=p(q_1+q_2)$.
Then $qq_i \leq d_ipq_i=0$, $p+q_i=p+d_i(p+d_j) \geq q_j$
whence $q+q_i =(q_1+q_2)(p+q_i) =q_1+q_2$. It follows
that $q/0$ is projective to $q_i/0$ whence $q_i \in I$. 
Then also $p_i= q_i+d_ip \in I$. 
Thus, in both cases  by the inductive hypothesis, there
exist  $u_{k1},...,u_{km} \in I$ with
$u_{ki} \leq  c_{i}$, $i = 1, ..., m$, and 
$p_k \leq f_{k}(u_{k1},..., u_{km})$ for $k=1,2$.  Set 
$u_{i} = u_{1i}+u_{2i}$, for $i = 1,..., m$.
and notice that $p \leq p_{1} + p_{2} \leq f(u_{1},...,u_{m})$.
The converse follows from monotonicity of lattice
polynomials. $\pfbox$

\xc{
\setlength{\unitlength}{0.00500000in}
\begingroup\makeatletter\ifx\SetFigFont\undefined%
\gdef\SetFigFont#1#2#3#4#5{%
  \reset@font\fontsize{#1}{#2pt}%
  \fontfamily{#3}\fontseries{#4}\fontshape{#5}%
  \selectfont}%
\fi\endgroup%
{\renewcommand{\dashlinestretch}{30}
\begin{picture}(893,695)(0,-10)
\put(120,385){\blacken\ellipse{10}{10}}
\put(120,385){\ellipse{10}{10}}
\put(200,465){\blacken\ellipse{10}{10}}
\put(200,465){\ellipse{10}{10}}
\put(280,545){\blacken\ellipse{10}{10}}
\put(280,545){\ellipse{10}{10}}
\put(360,625){\blacken\ellipse{10}{10}}
\put(360,625){\ellipse{10}{10}}
\put(440,545){\blacken\ellipse{10}{10}}
\put(440,545){\ellipse{10}{10}}
\put(520,465){\blacken\ellipse{10}{10}}
\put(520,465){\ellipse{10}{10}}
\put(600,385){\blacken\ellipse{10}{10}}
\put(600,385){\ellipse{10}{10}}
\put(520,305){\blacken\ellipse{10}{10}}
\put(520,305){\ellipse{10}{10}}
\put(440,385){\blacken\ellipse{10}{10}}
\put(440,385){\ellipse{10}{10}}
\put(360,465){\blacken\ellipse{10}{10}}
\put(360,465){\ellipse{10}{10}}
\put(280,385){\blacken\ellipse{10}{10}}
\put(280,385){\ellipse{10}{10}}
\put(200,305){\blacken\ellipse{10}{10}}
\put(200,305){\ellipse{10}{10}}
\put(280,305){\blacken\ellipse{10}{10}}
\put(280,305){\ellipse{10}{10}}
\put(360,385){\blacken\ellipse{10}{10}}
\put(360,385){\ellipse{10}{10}}
\put(440,305){\blacken\ellipse{10}{10}}
\put(440,305){\ellipse{10}{10}}
\put(360,305){\blacken\ellipse{10}{10}}
\put(360,305){\ellipse{10}{10}}
\put(280,225){\blacken\ellipse{10}{10}}
\put(280,225){\ellipse{10}{10}}
\put(280,145){\blacken\ellipse{10}{10}}
\put(280,145){\ellipse{10}{10}}
\put(360,145){\blacken\ellipse{10}{10}}
\put(360,145){\ellipse{10}{10}}
\put(340,185){\blacken\ellipse{10}{10}}
\put(340,185){\ellipse{10}{10}}
\put(440,145){\blacken\ellipse{10}{10}}
\put(440,145){\ellipse{10}{10}}
\put(440,225){\blacken\ellipse{10}{10}}
\put(440,225){\ellipse{10}{10}}
\put(440,145){\blacken\ellipse{10}{10}}
\put(440,145){\ellipse{10}{10}}
\put(360,65){\blacken\ellipse{10}{10}}
\put(360,65){\ellipse{10}{10}}
\put(360,225){\blacken\ellipse{10}{10}}
\put(360,225){\ellipse{10}{10}}
\put(361,24){\blacken\ellipse{10}{10}}
\put(361,24){\ellipse{10}{10}}
\put(200,145){\blacken\ellipse{10}{10}}
\put(200,145){\ellipse{10}{10}}
\put(520,145){\blacken\ellipse{10}{10}}
\put(520,145){\ellipse{10}{10}}
\put(200,145){\blacken\ellipse{14}{14}}
\put(200,145){\ellipse{14}{14}}
\put(200,305){\blacken\ellipse{14}{14}}
\put(200,305){\ellipse{14}{14}}
\put(120,385){\blacken\ellipse{14}{14}}
\put(120,385){\ellipse{14}{14}}
\put(600,385){\blacken\ellipse{14}{14}}
\put(600,385){\ellipse{14}{14}}
\put(520,305){\blacken\ellipse{14}{14}}
\put(520,305){\ellipse{14}{14}}
\put(440,305){\blacken\ellipse{14}{14}}
\put(440,305){\ellipse{14}{14}}
\put(280,305){\blacken\ellipse{14}{14}}
\put(280,305){\ellipse{14}{14}}
\put(340,185){\blacken\ellipse{14}{14}}
\put(340,185){\ellipse{14}{14}}
\put(520,145){\blacken\ellipse{14}{14}}
\put(520,145){\ellipse{14}{14}}
\put(360,625){\blacken\ellipse{14}{14}}
\put(360,625){\ellipse{14}{14}}
\thicklines
\put(494.819,125.336){\arc{338.154}{2.5136}{3.5086}}
\path(360,625)(120,385)
\put(542,289){\makebox(0,0)[lb]{\smash{{{\SetFigFont{12}{14.4}{\familydefault}{\mddefault}{\updefault}$p_2=d_2p\oplus q_2$}}}}}
\path(360,625)(600,385)
\path(120,385)(360,145)
\path(600,385)(360,145)
\path(280,545)(520,305)
\path(200,465)(440,225)
\path(440,545)(200,305)
\path(360,465)(360,305)
\path(280,225)(280,145)
\path(360,145)(360,65)
\path(440,225)(440,145)
\path(360,305)(360,225)
\path(360,225)(280,145)
\path(280,145)(360,65)
\path(440,145)(360,65)
\path(360,225)(440,145)
\path(360,385)(440,305)
\path(440,305)(360,225)
\path(280,305)(340,185)
\path(440,305)(340,185)
\path(360,65)(360,25)
\path(200,305)(200,145)
\path(280,305)(200,145)
\path(280,305)(520,145)
\path(520,305)(520,145)
\path(200,145)(360,25)
\path(520,145)(360,25)
\path(520,465)(280,225)
\put(162,114){\makebox(0,0)[lb]{\smash{{{\SetFigFont{12}{14.4}{\familydefault}{\mddefault}{\updefault}$q_1$}}}}}
\put(241,314){\makebox(0,0)[lb]{\smash{{{\SetFigFont{12}{14.4}{\familydefault}{\mddefault}{\updefault}$q_1+q_2$}}}}}
\put(454,296){\makebox(0,0)[lb]{\smash{{{\SetFigFont{12}{14.4}{\familydefault}{\mddefault}{\updefault}$p$}}}}}
\put(353,174){\makebox(0,0)[lb]{\smash{{{\SetFigFont{12}{14.4}{\familydefault}{\mddefault}{\updefault}$q$}}}}}
\put(527,119){\makebox(0,0)[lb]{\smash{{{\SetFigFont{12}{14.4}{\familydefault}{\mddefault}{\updefault}$q_2$}}}}}
\put(372,0){\makebox(0,0)[lb]{\smash{{{\SetFigFont{12}{14.4}{\familydefault}{\mddefault}{\updefault}$0$}}}}}
\put(0,290){\makebox(0,0)[lb]{\smash{{{\SetFigFont{12}{14.4}{\familydefault}{\mddefault}{\updefault}$d_1p\oplus q_1=p_1$}}}}}
\put(0,405){\makebox(0,0)[lb]{\smash{{{\SetFigFont{12}{14.4}{\familydefault}{\mddefault}{\updefault}$f_1(x_1,\ldots)=d_1$}}}}}
\put(155,645){\makebox(0,0)[lb]{\smash{{{\SetFigFont{12}{14.4}{\familydefault}{\mddefault}{\updefault}$f(x_1, \ldots)=f_1(x_1, \ldots)+f_2(x_1,\ldots)$}}}}}
\put(617,376){\makebox(0,0)[lb]{\smash{{{\SetFigFont{12}{14.4}{\familydefault}{\mddefault}{\updefault}$d_2=f_2(x_1, \ldots)$}}}}}
\end{picture}
}
}

\begin{corollary} Let $M$ be a complemented modular
 lattice with orthogonality $\perp$  and $I$ a neutral ideal such that for each
$a>0$ there is $p \in I$, $a \geq p>0$.
Then an orthoimplication holds in $M$ if and only if it holds
in all $[0,u]$, $u \in I$.
\lab{23}
\end{corollary}
{\em Proof.}
Consider an orthoimplication given  by $f$ which is not valid in $M$.
  There exist $x_{1} \perp y_{1}, ..., x_{n} \perp y_{n}$
so that $f(x_{1},y_{1},...,x_{n},y_{n})=a> 0$
whence $0<p\leq a$ with $p \in I$.  By    \ref{22},
there exist $u_{i},v_{i} \in I$ with $u_{i} \leq  x_{i}$, 
$v_{i} \leq  y_{i}$, for 
$i = 1, ..., n$ and $f(u_{1},v_{1}, ..., u_{n}, v_{n})\geq p > 0$.
But $u_{i} \leq  x_{i}$ 
and $v_{i} \leq  y_{i}$ give  $u_{i}
\perp v_{i}$, for $i = 1,...,n$.  Let
$u = \sum_{i=1}^{n} u_{i} + \sum_{i=1}^{n} v_{i}$. 
Then the orthoimplication  does not hold in $[0,u]$. $\pfbox$

\begin{corollary} Let $M$ be a algebraic lattice with orthogonality.
Then an orthoimplication holds in $M$ if and only if it holds
for all substitutions with compact elements.
\lab{23b}
\end{corollary}

\begin{lemma}
Within the variety of orthomodular lattices, each ortholattice identity
is equivalent to an orthoimplication in terms of the canonical orthogonality.
\lab{oimp1}
\end{lemma}

\no
{\em Proof.}
Considering an identity $g=h$ in the language of ortholattices
we may replace the constants $0,1$ by
$uu'$ resp. $u+u'$, $u$ a new variable. Also, we may assume that 
$g \leq h$ is valid in all ortholattices.
  Due to DeMorgan's Laws and $x'' = x$,
there is a lattice term $f(x_{1},y_{1}, ..., x_{n},y_{n})$
such that $hg'(x_{1},...,x_{n})
= f(x_{1},x_{1}',...,x_{n},x_{n}')$ holds in all ortholattices.
If $g = h$ holds in the orthomodular lattice $L$, and $x_{i} \perp y_{i}$,
$i = 1,...,n$, then $0 = f(x_{1},x_{1}',..., x_{n}, x_{n}')
\geq f(x_{1},y_{1},...,x_{n},y_{n})$ and the orthoimplication
\[  x_{1} \perp y_{1} \wedge ... \wedge x_{n} \perp y_{n}
\rightarrow f(x_{1},y_{1},...,x_{n},y_{n}) = 0\]
holds in $L$.  Conversely, if this orthoimplication
holds, then $f(x_{1},x_{1}',...,x_{n},x_{n}') = 0$ holds
in $L$ and, consequently, $g = h$ holds in $L$. $\pfbox$
With  \ref{23},     \ref{oimp1}, and \ref{neu}  one gets

\begin{corollary}
For a subdirectly irreducible MOL $L$ with minimal congruence $\mu$
the variety $V(L)$ is generated by the
simple  interval subalgebras $[0,u]$ of $L$, $u/0\in\mu$. \lab{molvar}
In particular, every variety of MOLs is generated by its simple
 members.
\end{corollary}

\begin{corollary}
The variety $V(L)$ of an atomic MOL $L$ is generated by
the interval subalgebras $[0,u]$, $u\in L_{fin}$.
\lab{atvar}
\end{corollary}
MOLs in the variety generated by atomic MOLs
(i.e. by finite dimensional MOLs) will be called
{\em proatomic}.

\begin{corollary} If $L \subseteq M$ is a  geometric representation 
of an MOL then   the orthoimplications of $M$ are
valid in $L$ and
 $L$ belongs to the variety generated by
$\hat{L}$ resp.   \lab{eqrep}
the set of interval subalgebras $[0,u]$ of $M$, $u \in M_{fin}$
\end{corollary}
{\em Proof}.
 Use Lemma \ref{HSP} and the fact that $L$ is a weak substructure of $M$.
Also, use 
\ref{23}, \ref{fin}, and
 \ref{oimp1}.  $\pfbox$

\subsection{Atomic extension}

For the proof of Thm \ref{atex} we need the following Lemma. The concept
of {\em neutral filter} is the dual of ``neutral ideal''.
We write $p\leq F$ if $p \leq x$  for all $x \in F$.

\begin{lemma} Let  $L,M$ be MOLs, $L$ a subalgebra of $M$, and 
 $F$ a neutral filter of $L$. Consider $a,b \in L$ and $p \in P_M$ such that
$ab=0$, $p \leq a+b$, and $p \leq F$. 
Then  \lab{neufi}
$a(p+b),\,b(p+a) \leq F$.
\end{lemma}
{\em Proof}. 
In view of restriction to interval subalgebras, we may assume $a+b=1$. 
Let $q=a(p+b)$ and $r=b(p+a)$  and $\theta$ the congruence
associated with $F$. Consider $x \in F$, i.e. $x\,\theta\,1$
and $p\leq x$.
 Let
\[ y=(a+xb)(b+x) \geq q,\;\; z=(b+xa)(a+x)\geq r\]
By modularity, $x,y,z$ coincide or are the atoms
of a sublattice of height $2$. In particular, all its quotients are
in $\theta$ 
whence $1/y \in \theta$ and $y \in F$. From $p \leq F$ it follows $p \leq y$ and
thus
$r \leq p+q \leq y$. Hence $r \leq yz \leq x$ and $q \leq x$, symmetrically.
$\pfbox$\\
{\em Proof  of  Thm}. \ref{atex}. 
(2) and (3) are equivalent by Cor.\ref{canrep} and Prop.\ref{geovar},
and imply (1) by Cor.\ref{atvar}.
The class of MOLs admitting an atomic extension 
 contains all finite dimensional ones  and is closed
under  subalgebras and direct products.
Thus, to prove that (1) implies (2) we have to show that this class is
closed under homomorphic images, too. Consider
a subalgebra  $L$ of an atomic MOL $M$ and congruence $\theta$ on $L$
with associated neutral filter $F$.  
Define
\[ Q=\{p\in P_M\mid p \leq F\},\;\; \eta: L/\theta \rightarrow L(Q),\;\; \eta(a/\theta)=\{p \in Q\mid p \leq a\}\]
Then $Q$ is a subgeometry of $P_M$ with polarity $\perp$, obviously, $\eta$ is meet
preserving and $\eta(a/\theta) \perp \eta(a'/\theta)$.
If $a/b \in \theta$ then $b =ac$ for some $c \in F$ whence
$a \geq p \in Q$ implies $p \leq b$; thus, $\eta$ is well defined. 
The proof that $\eta$ preserves joins follows Frink:  
Given $a,b \in L$ choose $\ti{b}$ such that $a+b=a+\ti{b}$ and
$a \ti{b}=0$. Consider $p \in \eta((a+b)/\theta)$, $p \not\in \eta(a/\theta)$ and
$p \not\in \eta(\ti{b}/\theta)$. Then by Lemma \ref{frink0}
$p, a(p +\ti{b})$, and $\ti{b}(p+a)$ are collinear elements of
$P_M$.  By Lemma \ref{neufi} they are in $Q$, whence
$p \in \eta(a/\theta) +\eta(\ti{b}/\theta) \subseteq \eta(a/\theta) +\eta(b/\theta)$.

Finally, consider $a/0 \not\in\theta$ which means $ac >0$ for all $c \in F$. Thus,
since $F$ is closed under finite meets, for any finite
$C \subseteq F$ we have $x \in M$ such that
$x\leq c$ for all $c \in C$. In other words, the
set 
\[ \Phi_a(x)=\{0<x \leq ac\mid c \in F\} \] 
 of formulas with parameters in $L$ is finitely satisfiable in $M$.
By the Compactness Theorem of First Order Logic,
$M$   has an elementary extension $M'$ such that
each $\Phi_a(x)$ is satisfiable  in $M'$, i.e. there is
$x \in M'$ with $0<x\leq ac$ for all $c \in F$.
Replacing $M$,  we may assume $M=M'$.
 Since $M$ is atomic, we get $p \in P_M$ with $p \leq x$ and
then $p \in Q$ by definition. 
Thus $\eta(a/\theta)>0$ 
which proves that $\eta$ is a geometric representation. $\pfbox$
With Cor.\ref{odec} we obtain

\begin{corollary} Every proatomic MOL  has a geometric representation
in an orthogonal union of spaces $P_i$, each of is given by a vector space $V_i$
over a $*$-division-ring $D_i$ with
anisotropic $*$-hermitian form $\Phi_i$
- or possibly  of height $3$ if $L$ is not Arguesian. Every  subdirectly
irreducible  proatomic MOL  has a representation with a single $P_i=P$.
\end{corollary}
Von Neumann \cite{NEU2} constructs a continuous, simple, atomless  MOL as
the  
metric completion of a direct union of finite dimensional MOLs.
Since the metric completion amounts to
a homomorphic image of a subalgebra of a direct
power, this MOL is proatomic. The finite dimensional MOLs are the 
 ${\cal L}(\Bbb{R}^{2^n} _{\Bbb{R}})$
and the union is formed with respect to the canonical  embedding maps
\[\phi_n:{\cal L}(\Bbb{R}^{2^n} _{\Bbb{R}})\rightarrow 
{\cal L}(\Bbb{R}^{2^{n+1}}_{\Bbb{R}}),\quad  \dim \phi_n x= 2\cdot \dim x\]

\subsection{Interpretation of $*$-ring identities}

Frames have played a crucial r\^{o}le in the equational theory of modular
lattices - due to the fact that the modular lattice
freely generated by an $n$-frame is projective 
with respect to onto-homomorphisms.
The analogous result holds according to 
 Mayet and Roddy \cite{MR} for  orthogonal $n$-frames  within the
variety of
relative  MOLs.
The following is the basis for connecting the equational
theories of MOLs and 
$*$-regular rings.

\begin{lemma}
There exist ortholattice-polynomials $t(x)$ and $x\st$
with constants   from $\f{a}$
such that for any MOL with spanning orthogonal $n$-frame ($n \geq 4$ or
Arguesian)
 \[(r^*)_{12}=(r_{12})\st,\quad 
\forall x.\;t(x)\in R_{12}, \quad t(r_{12})=r_{12} \]
\end{lemma}
{\em Proof.}
Indeed, $e_1-e_2r \perp e_1r^*\al_2+e_2$ whence
$(e_1r^*\al_2+e_2)R \leq r'_{12}(a_1+a_2) \in R_{21}$ 
and equality follows by modularity. Thus
 \[(-r^*\al_2)_{21}=r'_{12}(a_1+a_2),\;\;(-\al_2)_{21}=a'_{12}(a_1+a_2)\]
$t(x)$ is provided by  the  lattice term
\[l(x,x',z_1, \ldots ,z_n)=(x'+x(x'+\sum_{j \neq 2}z_j))(x+ x'\sum_{j \neq 2}z_j)(z_1+z_2)\]
Observe that $\hat{x}=(x'+x(xz_2)')(x+(x+z_2)')$ is a complement of $x$
in $[x(xz_2)',x+(x+z_2)']$ whereas $x(xz_2)'$ is a complement
of $xz_2$ in $[0,x]$ and $x+(x+z_2)'$ a complement of $x+z_2$ in $[x,1]$.
Therefore, $\hat{x}$ is a complement of $z_2$ and $\hat{x}(z_1+z_2)$
a complement of $z_2$ in $[0,z_1+z_2]$.
 Now, for any given  spanning orthogonal frame $\f{a}$
one has $l(x,a_1, \ldots ,a_n)=\hat{x}(a_1+a_2)$ and  it follows 
$l(x,x',a_1, \ldots,a_n) \in R_{12}$.  $\pfbox$
Combining this with the Mayet-Roddy  terms providing the orthogonal
frame and the polynomials yielding the structure of the  coordinate ring, one obtains
the following.

\begin{theorem} For every $*$-ring identity $\al$ there
is an MOL-identity $\hat{\al}$ such that for every $*$-ring $R$ 
associated with a $*$-regular matrix ring $R_3$, the identity
$\al$ holds in $R$ if and only if $\hat{\al}$ holds in
$\ov{L}(R_3)$.
\end{theorem}

\subsection{Generating frames}
 The
subdirectly  irreducible  frame generated objects  have been
determined for $n \geq 4$ (resp. Arguesian)  
modular lattices (\cite{Hfr}). The analogous 
task appears intractable for MOLs.
The starting point was the construction of a 
$3$-frame generated height $6$ MOL by B.M\"uller.

Let  $\f{E}^3$ denote the canonical $3$-frame of $L={\cal L}((R_2)^3_{R_2})$.
 The canonical isomorphism between 
 ${\cal L}(R^2_R)$ and ${\cal L}(R_{2R_2})$  gives rise to an 
isomorphism of ${\cal L}(R^6_R)$ onto 
$L$ mapping the canonical $6$-frame $\f{E}^6$ onto
$\tilde{\f{E}}^6:
  \tilde{e}^6_iR_2, (\tilde{e}^6_i-\tilde{e}^6_j)R_2$
where
\[\tilde{e}^6_i=e^3_i \pmatrix{ 1&0\cr 0&0 },\;
\tilde{e}^6_{i+3}=e^3_i\pmatrix{ 0&0\cr 1&0 } \mbox{ for } i=1,2,3\]
Let  $\Q$ be the field of rational numbers.
\begin{lemma}
 Let $R$ be a finite dimensional $\Q$-algebra
 and $a,b$  invertible elements of 
$R$ such that all $a+(1-\frac{1}{2^k})b$ are invertible. Let $S$
be
 generated by $a,b$ under  ring operations and inversion (as far as inverses exist)  and $M$ be 
sublattice  of $L$ generated by $\f{E}^3$ and $A_{12}, B_{13}$
where \lab{M}
\[ A= \pmatrix{a+b&b\cr b&2b},\;\; B=\pmatrix{ b&b\cr b&2b} \]
Then $A$,$B$, and all
$A+(1-\frac{1}{2^k})B$  are invertible in $R_2$, $\tilde{\f{E}}^6 \subseteq M$  and $C_{ij} \in M$
for every matrix $C \in S_2$. 
\end{lemma}
{\em Proof.}
\[\pmatrix{1 &-\frac{1}{2}\cr 0 &1} (A+(1-\frac{1}{2^k})B)
 =\pmatrix{a +(1-\frac{1}{2^{k+1}})b& 0\cr (2-\frac{1}{2^{k}})b  &2(2-\frac{1}{2^{k}})b},\;\;
\pmatrix{ 1 & 0\cr -1 &1}
B= \pmatrix{b & b \cr 0 & b}\]
 Calculating in $M$ resp. $R(M,\f{E}^3)$  
we get
\[\pmatrix{a&0\cr0&0}=A-B, \quad \tilde{E}^6_2=
E^3_2 \cap (E^3_1+\pmatrix{a&0\cr0&0}_{12}), \quad
\tilde{E}^6_4= \pmatrix{a&0\cr0&0}_{12}\cap E^3_1  \]
whence $\tilde{E}^6_i=E_i\cap (\tilde{E}^6_2+E_{2i})$ and
$\tilde{E}^6_{ij}=E_{ij}\cap (\tilde{E}^6_i+\tilde{E}^6_j)$ for $i,j\leq 3$ and,
similarly, for $i,j \geq 4$.
In particular 
\[ \pmatrix{1&0\cr0&0}_{12} = (\tilde{E}^6_{12}+\ti{E}^6_4)\cap (E^3_1+E^3_2) \in  R(M,\f{E}^3)_{12}\]
Thus we have in  $R(M,\f{E}^3)$
\[ \pmatrix{b&0\cr 0&0} =\pmatrix{1&0\cr0&0}B\pmatrix{1&0\cr0&0},\;\;\mbox{ and }
 \pmatrix{c&0\cr 0&0}  \mbox{ for all } c\in S\]
since 
\[ \pmatrix{c^{-1}&0\cr 0&0}_{12} = (\pmatrix{c&0\cr 0&0}_{21}+\ti{E}^6_4
+\ti{E}^6_5)\cap ({E}^3_1+\ti{E}^6_2) \]
Now, we get in $R(M,\f{E}^3)$
\[\pmatrix{0&b\cr b&2b}= B-\pmatrix{a+b&0\cr0&0}\]\[\pmatrix{0&0\cr1&0} =\pmatrix{0&b\cr b&2b}\pmatrix{b^{-1}&0\cr 0&0},\;\;
\pmatrix{0&1\cr0&0} =\pmatrix{b^{-1}&0\cr 0&0} \pmatrix{0&b\cr b&2b}\]
whence all of $S_2$. Moreover, we have $\tilde{\f{E}}^6 \subseteq M$ from
\[ \tilde{E}^6_{15} 
= (\pmatrix{0&0\cr 1&0}_{12}+\ti{E}^6_4) \cap (\ti{E}^6_1+\ti{E}^6_5) \in M. \quad\quad \pfbox \]   
Define $R(1)= \Q$, $A_1=B_1=(1)$
and, inductively,
\[R(k+1)=R(k)_2,\;\;A_{k+1}=\pmatrix{A_{k}+B_k&B_k\cr
B_k&2B_k},\; 
B_{k+1}=\pmatrix{B_{k}&B_k\cr B_k&2B_k}.\]
\begin{lemma} The lattice ${\cal L}(R(k)^3_{R(k)})$ is generated by its canonical
$3$-frame together with $(A_k)_{12}$ and $(B_k)_{13}$. \lab{M1} 
\end{lemma}
{\em Proof.}
The case $k=1$ is well known, cf \cite{B67}.  
Now, in the inductive step $k \to k+1$ 
we use \ref{M} with $R=R(k)$, $a=A_k$, $b =B_k$.
We have ${\cal L}(R^3_R)$ embedded via $\phi$ into $L$ with the canonical $3$-frame
 mapped onto $\ti{E}^6_i,\ti{E}^6_{ij}$, $i,j \leq 3$, all
contained in $M$. Also,
\[ \phi a_{12}=(\ti{E}^6_1+\ti{E}^6_2)\cap(\pmatrix{a&0\cr 0&0}_{12}+\ti{E}^6_4),\;\;
\phi b_{13}=(\ti{E}^6_1+\ti{E}^6_3)\cap (\pmatrix{b&0\cr 0&0}_{13}+\ti{E}^6_4)\]
belong to $M$.
By the inductive hypothesis,  ${\cal L}(R^3_R)$
is generated by $a_{12}, b_{13}$ together with the canonical  $3$-frame.
Thus, all of the image belongs to $M$ and so does
\[ \pmatrix{c &0\cr0 &0}_{12}= (\phi c_{12} +\ti{E}^6_4) \cap (E^3_1+E^3_2),\;
\mbox{ where } c \in R\]
As above, we get $R(M,\f{E}^3)=R_2 =R(k+1)$ and it follows
$M=L$. $\pfbox$

\begin{proposition} \lab{M2}
For all $n=2^k$ 
there is a positive definite symmetric form on the vector space
$\Q^{3n}$
 such that the image  of the canonical frame 
of ${\cal L}((\Q_n)^3_{\Q_n})$ is a  generating orthogonal $3$-frame in
the MOL ${\cal L}(\Q^{3n}_{\Q})$. 
\end{proposition}
{\em Proof}. Start with $a,b>0$ in $\Q$ and
consider the above defined  $A_k$, $B_k$ as 
$2^k\times 2^k$-matrices over $\Q$. 
Induction and the congruence transformations
\[\pmatrix{1 &-\frac{1}{2}\cr 0 &1} A \pmatrix{1& 0\cr -\frac{1}{2} &1} =\pmatrix{a +\frac{1}{2} b& 0\cr 0 &2b},\;\;
\pmatrix{ 1 & 0\cr -1 &1}
B\pmatrix{ 1 & -1\cr 0 &1}= \pmatrix{b &0 \cr 0 & b}\]
 show that both are positive definite symmetric matrices.
 Endow 
 $\Q^{3n}$ with the form
given by the positive definite block matrix
\[ \pmatrix{I_k &O&O\cr O&A^{-1}_k&O \cr O&O& B^{-1}_k} \]
and ${\cal L}(\Q^{3n})$ with the induced orthocomplementation. 
Under the isomorphism $\psi: {\cal L}((\Q_n)^3_{\Q_n}) \imp {\cal L}(\Q^{3n}_{\Q})$  the image
 $\ti{\f{E}}^3$ of the canonical $3$-frame $\f{E}^3$ consists of
\[  \ti{E}^3_i =\sum_{t=1}^{n} e_{t+(i-1)n}\Q,\;\;
\ti{E}^3_{ij} =\sum_{t=1}^{n} (e_{t+(i-1)n}- e_{t+(j-1)n})\Q .\]
The $\ti{E}^3_i$, $i=1,2,3$ are pairwise orthogonal.
Moreover, in $R({\cal L}(\Q^{3n}),\ti{\f{E}}^3)$ we have
\[ \psi (A_k)_{12}=\ominus_{12}((\ti{E}^3_{12})'\cap (\ti{E}^3_1+\ti{E}^3_2)), \quad 
\psi (B_k)_{13}=\ominus_{13}((\ti{E}^3_{13})'\cap  (\ti{E}^3_1+\ti{E}^3_3)).\]
By \ref{M1} ${\cal L}((\Q_n)^3_{\Q_n})$ is generated as a lattice by
$\f{E}^3$ and $(A_k)_{12}, (B_k)_{13}$.  
Hence, the MOL ${\cal L}(\Q^{3n}_{\Q})$ is generated by $\ti{\f{E}}^3$.
$\pfbox$
Let $L$ a non-principal ultraproduct of the $L_k$, $k\geq 1$, and
let $\f{a}$ correspond to the $\f{a}_k$. Then in the sublattice
generated by $\f{a}$, for any $k$ one has $x_1> \ldots >x_k$ with 
 the $x_i/x_{i+1}$ pairwise projective. Hence
\begin{corollary}
There is a subdirectly irreducible MOL of infinite height
generated by an orthogonal $3$-frame.  
\end{corollary}

\subsection{Word problems}
Finally, we consider quasi-identities
$\bigwedge_i s_i=t_i\imp s=t$ resp. their
model classes, called quasi-varieties.
Recall that there is a MOL \cite{BR} not in the quasi-variety generated by finite
dimensional MOLs.
The {\em  word problem} for a quasi-variety requires an
algorithm dealing with all finite presentations,
 i.e. a decision procedure for quasi-identities.

\begin{proposition}
Let {\cal Q} be any quasi-variety  of modular (ortho)lattices containing all
${\cal L}(\Q^n_{\Q})$, $n\geq 4$ (with orthogonality
given by the identity matrix).  
 Then {\cal Q} 
has unsolvable  word problem. 
\end{proposition}
{\em Proof.}
Let $\Lambda$ denote the
 set of all quasi-identities in the language of
semigroups.
According to Gurevich and Lewis \cite{GL} there is no recursive $\Ga \subseteq \Lambda$ such
that $\Ga$ contains all $\phi$ valid in all semigroups but
  none falsified in some finite semigroup. 
Associate with each  $\phi$ in $\Lambda$ 
a lattice quasi-identity $\hat{\phi}$ expressing that the semigroup variables
correspond to elements of the coordinate ring of a
 $4$-frame and translating semigroup relations
into lattice relations (cf \cite{Hemvd} for a similar translation).   
Here, $4$-frames with a family of elements of
the coordinate ring  have to be considered as systems of lattice generators
and relations (as defined by von Neumann \cite{NEU}).
 By the Coordinatization Theorem, 
the coordinate ring is indeed a ring under the intended operations.
Therefore, if  $\Ga$ is the set of all $\phi$ with $\hat{\phi}$ valid
in  {\cal Q} then  $\Ga$ contains all $\phi$ valid in all semigroups.
On the other hand, if $\phi$ is falsified in the finite
semigroup $S$, we represent $S$ as a 
subsemigroup of some matrix ring $\Q_n$, i.e. a subsemigroup of the
  coordinate ring of  ${\cal L}((\Q_n)^4_{\Q_n})$
with canonical $4$-frame. The lattice may be turned into an MOL
transferring the canonical orthogonality of ${\cal L}(\Q^{4n}_{\Q})$
via an isomorphism. 
Thus, $\hat{\phi}$ is falsified in $L$ which means $\phi \not\in \Ga$.
Now, assuming that {\cal Q} has decidable
quasi-identities  would yield that $\Ga$ is recursive, a contradiction. 
$\pfbox$ In the case of modular lattices, $\Q$  
can be replaced by any prime field.
The task of finding a particular finite presentation with 
unsolvable word problem is  substantially  more demanding. 
 It has been completed  for modular lattices with $5$ generators
by Hutchinson \cite{Hut} under the same assumption, for MOLs
with $3$ generators in \cite{R:89}  
for each quasi-variety  containing all
subdirectly irreducible MOLs of height 14.

\section{Discussion}

The Frink space of an MOL, $L$, is endowed with a canonical
anisotropic orthogonality $\perp$ satisfying all orthoimplications of $L$ 
according to Lemma \ref{HSP}. So, if $\perp$ is a polarity, Prop.\ref{geovar}
 provides 
a canonical atomic extension within the variety of $L$.
Unfortunately, the direct union of MOLs in the von Neumann
example constitutes a counterexample, already.

\begin{problem}
Characterize the MOLs for which the Frink embedding
provides a geometric representation.
\end{problem}
\begin{problem} Does every MOL admit an atomic extension? \lab{MOL}
\end{problem}
\begin{problem} Does every proatomic MOL admit an extension within
its variety? \lab{varex}
\end{problem}
For the last problem, the following concept might be helpful.
Call 
a  geometric representation  $\eta:L \imp M= {\cal S}(P)$    {\em orthogonally
separating} 
if  for all $u\perp v$ in the same component of $M_{fin}$ there is  some $a \in L$ 
such that
$u \leq \eta(a) $ and $v \leq \eta(a')$. 
Observe that 
then $L$ and ${\cal S}(L)$ satisfy the same orthoimplications and
Prop.\ref{geovar} provides an atomic extension in $V(L)$.
Also, the class of MOLs admitting such a representation is closed
under subdirect products. Thus, 
a positive answer to the following problem would also 
imply that for \ref{varex}.
\begin{problem} Is the class of MOLs admitting
an orthogonally separating representation closed under
homomorphic images?
\end{problem}
Actually, the original motivation for this research was 
the following question  partly answered in Roddy \cite{MR:87}. 
\begin{problem} Which MOL varieties, not generated by an MO$_\kappa$,
do contain a projective plane?
\end{problem}
G.Bruns \cite{GB:83} conjectured that it is true
for all varieties.  But the answer for proatomic
varieties is open  as well. Results of \cite{R:89} suggest   
that the equational theory of MOLs with suitable
 bound on the height of irreducible factors should be  undecidable. 
\begin{problem} Is the equational theory
of (proatomic) Arguesian  MOLs decidable?
\end{problem}

\abs
\begin{conjecture}  
The von Neumann example of a continuous geometry admits a geometric representation over
an elementary extension of the reals.
\end{conjecture}
As we have seen, the  von Neumann example is proatomic.
How far does this extend to abstract continuous
geometries - a positive answer 
could be seen as a kind of  construction for these.
Recall, that by Kaplansky \cite{KAP} and Amemiya and Halperin
\cite{AH} every countably complete MOL is continuous
and every continuous MOL is `finite'.
\begin{problem} Is every `finite'  (continuous, countably complete, complete)
MOL pro\-atomic? Do such even  belong to the quasivariety
generated by finite dimensional MOLs?
\end{problem}
Recall, that the quasivariety generated by a class consists of the
 subalgebras of products of ultraproducts.
 In a quasivariety generated by modular lattices 
of finite height, no quotient may be projective to a proper
subquotient - a property shared with modular lattices admitting
a dimension function.     
Wehrung \cite{Wehr2} calls a lattice {\em normal} 
if projective $a,b$ with $ab=0$ are perspective.
Bruns and Roddy \cite{BR} provide 
an atomic  MOL which is not normal. 
\begin{problem} Is normality inherited by sub-MOLs?
\end{problem}
Also, the representing space is of interest. 
According to Gross \cite{GRO} p.65 every  
hermitian vector  space of countable dimension  admits an orthogonal basis.
\begin{problem} Can every Arguesian proatomic  MOL be
represented by means of spaces  having orthogonal bases?
\end{problem}

\abs
Concerning coordinatization,
one has to ask how far J\'{o}nsson's results \cite {J60}
for
complemented modular lattices extend to MOLs.
J\'{o}nsson  constructed  an example of a simple
coordinatizable lattice with no spanning $n$-frame ($n \geq 3)$ which 
lead him to consider `large partial $n$-frames', $n \geq 3$.
He showed that every complemented modular lattice $L$
with such frame  ($n \geq 4$ or  $L$ Arguesian)
  is   coordinatizable and that every simple $L$
of height $\geq 4$ contains such a frame. 
We suggest the following definition of an {\em orthogonal large partial
$n$-frame}: For given  $m\geq n\geq 3$
it is constituted by  orthogonal elements $a_i\,(1\leq i \leq m)$ such that
\[ \sum_{i=1}^ma_i=1,\quad  a_i \sim a_1 \mbox{ for } i \leq n,\quad  a_i \sim y_i \leq a_1
\mbox{ for } n<i\leq m\]
\begin{conjecture}
The analogues of J\'{o}nsson's results hold for MOLs with
orthogonal large partial $n$-frames.
\end{conjecture}
This would imply that every MOL-variety is generated 
by members of height $\leq 3$ and  members of the form $\ov{L}(R)$ with
simple $R$. Yet, even so  one might  fail to characterize 
coordinatizability.
\begin{conjecture} There are  subdirectly irreducible  coordinatizable
MOLs of height $\geq 3$ not containing an orthogonal large partial $n$-frame.
\end{conjecture}
\begin{problem} Can every Arguesian MOL be
embedded into the interval $[0,a_1]$ of an MOL
with orthogonal large partial $3$-frame?
\end{problem}
\begin{problem} Does every $*$-regular ring belong to the
$*$-ring variety generated by $*$-rings $R$ associated with 
$*$-regular  matrix rings
$R_3$?
\end{problem}

\abs 
For  $*$-regular rings, the following
 concept appears to reflect geometric representation of MOLs.
A {\em representation} of $*$-regular ring $R$
is given by a vector space $V_D$,
 a ring embedding  $\iota:R \imp \mbox{\rm End}(V_D)$, and a  $*$-hermitian form $\Phi$ on $V_D$
such that 
$\iota(r^*)$ is the  adjoint of $\iota(r)$ for all $r \in R$
The following is due to Kaplansky (cf  \cite{HE})
\begin{theorem} 
Primitive    $*$-regular rings  
with  minimal left ideal       are  representable.
\end{theorem}
Characterizing representability in terms of 
proatomic MOLs could provide a construction of representable
rings from artinian $*$-regular rings and shed light on 
the type $I_n$ and $II_1$ factors of von Neumann algebras.
\begin{conjecture}  Every subdirectly irreducible 
representable $*$-regular  ring 
can be embedded into a homomorphic image of a regular 
$*$-subring of an ultraproduct of artinian $*$-regular rings.
\end{conjecture}
\begin{problem} Is every  $*$-regular ring representable?
\end{problem}
The following two concepts are quite important in the
theory of regular rings: A ring is {\em unit regular} if
for every $a$ there is a unit $u$ such that $aua=a$.
A ring is {\em directly finite} if $xy=1$ always implies
$yx=1$.  Observe that every artinian regular ring is unit regular
and every unit regular ring  is directly finite. 
Moreover, a regular ring $R$ with  $n$-frame ($n \geq 2)$
in $\ov{L}(R)$ is unit regular if and only if
perspectivity is transitive in this lattice. 
The following is due to Handelman (see \cite{G})
\begin{problem}  Is every $*$ regular ring directly
finite or even unit regular?
\end{problem}
\begin{conjecture} 
If $R$ is $*$-regular and $\ov{L}(R)$ proatomic
with orthogonal large partial $n$-frame  then $R$ is unit regular.
Every representable ring is directly finite
- and unit regular, if simple. 
\end{conjecture}
Some positive evidence is given by the
following results of  Ara and Menal \cite{AM} and of Kaplansky \cite{KAP}
and Amemiya and Halperin \cite{AH}.
\begin{theorem} If $R$ is $*$-regular, then $xx^*=1$ implies $x^*x=1$.
If, in addition,
 $\ov{L}(R)$ is 
$\aleph_0$-complete then $R$ is unit regular.\end{theorem}

\end{document}